# NAVIGATION ON A POISSON POINT PROCESS[1]

By Charles Bordenave

*University of California*

On a locally finite point set, a *navigation* defines a path through the point set from one point to another. The set of paths leading to a given point defines a tree known as the *navigation tree*. In this article, we analyze the properties of the navigation tree when the point set is a Poisson point process on $\mathbb{R}^d$. We examine the local weak convergence of the navigation tree, the asymptotic average of a functional along a path, the shape of the navigation tree and its topological ends. We illustrate our work in the small-world graphs where new results are established.

## 1. Introduction.

1.1. *Navigation*: *definition and perspective.* Let $N$ be a locally finite point set and $O$ a point in $\mathbb{R}^d$, taken as the origin. For $X, Y \in \mathbb{R}^d$, $|X|$ will denote the Euclidean norm and $\langle X, Y \rangle$ the usual scalar product.

DEFINITION 1.1.  Assuming that $O \in N$, a *navigation* (with root $O$) is a mapping $\mathcal{A}$ from $N$ to $N$ such that $\mathcal{A}(O) = O$ and for all $X$ in $N$, there exists $k \in \mathbb{N}$ such that $\mathcal{A}^k(X) = O$.

A navigation is the ancestor mapping of a rooted tree, the *navigation tree* $\mathcal{T}_O = (N, E_O)$, which is defined by

$$(X, Y) \in E_O \qquad \text{if } \mathcal{A}(X) = Y \text{ or } \mathcal{A}(Y) = X.$$

In this paper, we examine decentralized navigation algorithms over a random point set. The concept of decentralization is loosely defined in a mathematical framework and we will introduce another property of these mappings,

Received February 2006; revised August 2007.
[1]Supported by the Science Foundation Ireland Research Grant SFI 04/RP1/I512.

*regeneration* (Definition 1.2). Roughly speaking, if the point set is a sample of a Poisson point process, then decentralization and regeneration coincide.

Navigation algorithms have emerged in different classes of problems. A first class is the small-world phenomenon. Kleinberg [16] showed that this phenomenon relies on the existence of shortcuts in a decentralized navigation on a geometric graph (for extension and refinements see Franceschetti and Meester [10], Ganesh et al. [8, 12]). A second field of application is computational geometry. Kranakis, Singh and Urrutia [18] have introduced compass routing (see Morin [22] for a review of this class of problems). The ideas of computational geometry may benefit the design of real-world networks: sensor and ad hoc networks (see Ko and Vaidya [17]) or self-organized overlay and peer-to-peer networks (see Plaxton, Rajaraman and Richa [27], Liebeherr, Nahas and Si [20] or Ganesh, Kermarrec and Massoulié [11]). Finally, in the probabilistic literature, some authors have examined decentralized navigation algorithms (under other names): the Markov path on the Delaunay graph in Baccelli, Tchoumatchenko and Zuyev [5], the Poisson forest of Ferrari, Landim and Thorisson [9] and the directed spanning forest introduced by Gangopadhyay, Roy and Sarkar [13] (see also Penrose and Wade [23] and Baccelli and Bordenave [3]). The aim of the present work is to find a unified approach to these problems. The mathematical material used in this work is a natural extension of the ideas developed in [3].

We start by giving three examples of navigation. Among these three, only the last will draw our attention. These examples are nevertheless useful to understand the context better. Let $\mathcal{G}$ be a connected graph $\mathcal{G} = (N, E)$ with $O \in N$. Our first example is the *shortest path* navigation. Assume that the edges are weighted and define $\pi(X) = \{X, X_1, \ldots, O\}$ as the connected path from $X$ to $O$ which minimizes the sum of the weights of the edges on the path [assuming that $\pi(X)$ is uniquely defined]. Then, $\mathcal{A}(X) = X_1$ is the shortest path navigation. In [28, 29], Vahidi-Asl and Wierman have studied the shortest path on the Delaunay triangulation of a Poisson point process; see also Pimentel [26]. On the complete graph of a Poisson point process and weight on the edge $(X, Y)$ taken as $|X - Y|^\beta, \beta > 2$, an in-depth analysis was performed by Howard and Newman [15]. The shortest path is the continuum analog of the *first passage percolation* on the regular $\mathbb{Z}^d$-lattice. The shortest path navigation has poor decentralization properties.

An irreducible recurrent *random walk* on $\mathcal{G}$ defines almost surely a decentralized navigation. The length of the path from $X$ to $O$ is the hitting time of $O$ starting from $X$. However, on an infinite graph, one might expect that the walk is null recurrent and that this navigation will not be efficient.

The *maximal progress navigation* on $\mathcal{G}$ is a greedy algorithm: $\mathcal{A}(X) = Y$ if $Y$ is the neighbor of $X$ closest to $O$. Note that this navigation is not a proper navigation on all graphs since it may fail to converge to $O$. If $\mathcal{A}$ is a navigation, the *progress* at $X$ is defined as
$$P(X) = |X| - |\mathcal{A}(X)|.$$

Throughout this paper, to simplify the analysis, we will restrict ourselves to navigation algorithms with a nonnegative progress for all $X \in N$.

1.2. *Directed navigation.* Let $e_1 \in S^{d-1}$. A directed navigation with direction $e_1$ is a mapping $\mathcal{A}_d$ from $N$ to $N$ such that for all $X$ in $N$, $\lim_{k \to \infty} \langle \mathcal{A}_d^k(X), e_1 \rangle = +\infty$. Note that if $\langle e_1, e_2 \rangle > 0$, then a directed navigation with direction $e_1$ may also be directed with direction $e_2$. In the sequel, a navigation with direction $e_1$, will be denoted by $\mathcal{A}_{e_1}$ in order to stress the direction of interest. The *directed progress* (with direction $e_1$) at $X$ is defined by
$$P_{e_1}(X) = \langle \mathcal{A}_{e_1}(X), e_1 \rangle - \langle X, e_1 \rangle.$$

A few examples of directed navigation may be found in the literature: the directed path on the Delaunay tessellation [5], the Poisson forest [9], the directed spanning forest ([13, 23]). On a graph $\mathcal{G}$, we also define the *maximal directed progress navigation* as the navigation which maximizes the directed progress.

A directed navigation is the ancestor mapping in a directed forest. The *directed navigation forest*, $\mathcal{T}_{e_1} = (N, E_{e_1})$, is defined by
$$(X, Y) \in E_{e_1} \qquad \text{if } \mathcal{A}_{e_1}(X) = Y.$$

1.3. *Regenerative navigation on a Poisson point process.* Throughout this paper, $N = \sum_{n \in \mathbb{N}} \delta_{T_n}$ is a *Poisson point process* (PPP) of intensity 1 on $\mathbb{R}^d$. If $\delta_X$ is the Dirac mass at $X \in \mathbb{R}^d$, we will define $N^O = N + \delta_O$ and $N^{O,X} = N + \delta_X + \delta_O$. From Slyvniak's theorem, $N^O$ (resp. $N^{O,X}$) is a PPP in its Palm version at $O$ [resp. $(O, X)$]. Up to a scaling, it is not restrictive to assume that the intensity of $N$ is 1.

In this paper, we analyze the properties of decentralized navigation trees on $N^O$. The decentralization property on a PPP is captured by the notion of regeneration. Regenerative navigation is defined via memoryless navigation. For a directed navigation $\mathcal{A}_{e_1}$, let $X_k = \mathcal{A}_{e_1}^k(X)$ and $\mathcal{F}_k = \sigma\{X_0, \ldots, X_k\}$ be the filtration associated with the process $(X_k)_{k \in \mathbb{N}}$. $\mathcal{A}_{e_1}$ is *memoryless* if, for all $X \in \mathbb{R}^d$, $k \in \mathbb{N}$,

(1) $$\mathbb{P}(X_{k+1} - X_k \in \cdot | \mathcal{F}_k) = \mathbb{P}(\mathcal{A}_{e_1}(O) \in \cdot).$$

In other words, a directed navigation is memoryless if the process $(X_{k+1} - X_k)_{k \in \mathbb{N}}$ is an i.i.d. sequence which does not depend on the initial position $X_0 = X$. In dimension 1, a memoryless directed navigation with nonnegative directed progress defines a renewal point process. Similarly, for a navigation $\mathcal{A}$ on $N^O$, let $X_k = \mathcal{A}^k(X)$. $\mathcal{A}$ is a *memoryless* navigation if, for all $X \in \mathbb{R}^d$, $k \in \mathbb{N}$,

(2) $$\mathbb{P}(X_{k+1} \in \cdot | \mathcal{F}_k) = \mathbb{P}(X_{k+1} \in \cdot | X_k).$$

the sequence $(X_k)_{k \in \mathbb{N}}$ is then a Markov chain with $O$ as absorbing state.

DEFINITION 1.2. A directed navigation $\mathcal{A}_{e_1}$ is *regenerative* if there exists a stopping time $\theta > 0$ such that $X \mapsto \mathcal{A}_{e_1}^\theta(X)$ is a memoryless directed navigation and the distribution of $\theta$ is independent of $X$.

A navigation $\mathcal{A}$ is *regenerative* if there exist $x_0 > 0$ and a stopping time $\theta > 0$ such that $X \mapsto \mathcal{A}^\theta(X)$ is a memoryless navigation and the distribution of $\theta$ is independent of $X$ for $|X| \geq x_0$.

The stopping time $\theta$ will be called a *regenerative time*. As we will see, the analysis of regenerative navigation algorithms is easy. Then, using coupling techniques, we will exhibit a method to prove that a decentralized navigation on a PPP is regenerative. This provides a tool to analyze decentralized navigation algorithms on a PPP. This is the cornerstone of this paper.

1.4. *Notation.* $B(X, r)$ is the open ball of radius $r$ and center $X$, and $S^{d-1} = \{X \in \mathbb{R}^d : |X| = 1\}$ is the $d$-dimensional hypersphere. $\overline{A}$ will denote the closure of the set $A \subset \mathbb{R}^d$ and $A^c = \mathbb{R}^d \setminus A$. The $d$-dimensional volume of $B(O, 1)$ is $\pi_d = \pi^{d/2}/\Gamma(d/2 + 1)$ and $\omega_{d-1} = 2\pi^{d/2}/\Gamma(d/2)$ is the $(d-1)$-dimensional area of $S^{d-1}$. For $e_1 \in S^{d-1}$ and $X \in \mathbb{R}^d$, $\mathcal{H}_{e_1}(X) = \{Y : \langle Y - X, e_1 \rangle > 0\}$. $\ell_0(\mathbb{R}^d)$ [resp. $\ell_\infty(\mathbb{R}^d)$] will denote the set of measurable functions $s : \mathbb{R}^d \to \mathbb{R}_+$ such that $\lim_{|X| \to \infty} s(X) = 0$ (resp. $+\infty$). For a real random variable $Z$, if $F(t) = \mathbb{P}(Z \leq t)$, then $\overline{F}(t) = \mathbb{P}(Z > t)$. We will write $Z_1 \leq_{st} Z_2$ if, for all $t \in \mathbb{R}$, $\mathbb{P}(Z_1 > t) \leq \mathbb{P}(Z_2 > t)$ (stochastic domination). For a discrete set $A$, $|A|$ is the cardinality of $A$.

We will use $C_0$ to denote a positive constant to be thought of as small and $C_1$ to denote a positive constant to be thought of as large. The exact values of $C_0$ and $C_1$ may change from one line to the next. However, $C_0$ and $C_1$ will never depend on parameters of the problem.

Let $X_k = \mathcal{A}^k(X)$. Throughout this work, we will pay attention to the path from $X$ to $O$: $\pi(X) = \{X_0, \ldots, X_{H(X)} = O\}$ and to the distance of $X$ from $O$ in the navigation tree,

$$H(X) = \inf\{k \geq 0 : \mathcal{A}^k(X) = O\}.$$

1.5. *Examples.* For a few other examples, see [7].

*Small-world navigation.* The *small-world graph* is a graph $\mathcal{G} = (N^O, E)$ such that for all $X, Y \in N^O$, $(X, Y) \in E$ with probability $f(|X - Y|)$, independently of everything else, where $f$ is a nonincreasing function with values in $[0, 1]$. We will assume, as $t$ tends to infinity, that

$$f(t) \sim ct^{-\beta},$$

with $c > 0$ and $\beta > 0$. More formally, to each pair $(X, Y) \in N^O \times N^O$, $X \neq Y$, we associate an independent random variable $U(X, Y)$, uniform on $[0, 1]$, independent of $N$ and such that

$$(X, Y) \in E \quad \text{if } U(X, Y) \leq f(|X - Y|).$$

It is easily checked that the degree of a vertex $X$ is a.s. finite if and only if $\beta > d$. The small world graph $\mathcal{G}$ is sometimes referred to as the *long range percolation graph*. *Small-world navigation* is the maximal progress navigation on $\mathcal{G}$:

$$\mathcal{A}(X) = \arg\min\{|Y| : (X, Y) \in E\}.$$

As such, the small-world graph has isolated points and navigation is ill-defined on nonconnected graphs. To circumvent this difficulty, given $N$, for each $X \in N$, we condition the variables $(U(X, Y), Y \in N^O \cap B(O, |X|))$, on the event that $V_X = \{Y \in N^O \cap B(O, |X|) : (X, Y) \in E\}$ is not empty. This conditioning on $U$ is not problematic since, given $N$, for all $X \neq X'$ in $N$, the variables $(U(X, Y), Y \in N^O \cap B(O, |X|))$ and $(U(X', Y), Y \in N^O \cap B(O, |X'|))$ are independent. In particular, the sets $V_X$ and $V_{X'}$ are independent given $N$.

If $\beta > d$, the directed navigation with direction $e_1$ is defined similarly:

$$\mathcal{A}_{e_1}(X) = \arg\max\{\langle Y, e_1 \rangle : (X, Y) \in E\}.$$

The directed navigation is properly defined if the set of neighbors of $X$ in $\mathcal{H}_{e_1}$ is a.s. nonempty and finite. Hence, when dealing with $\mathcal{A}_{e_1}$, we will assume that the variables $U$ are conditioned on the event that a positive directed progress is feasible at any point $X$ of $N$.

*Compass routing on the Delaunay triangulation.* *Compass routing* was introduced by Kranakis et al. in [18]; see also Morin [22]. Let $\mathcal{G} = (N^O, E)$ denote a locally finite connected graph. Compass routing on $\mathcal{G}$ to $O$ is a navigation defined by

$$\mathcal{A}(X) = \arg\max\left\{\left\langle \frac{X}{|X|}, \frac{X - Y}{|X - Y|} \right\rangle : (X, Y) \in E\right\}.$$

$\mathcal{A}(X)$ is the neighboring point of $X$ in $\mathcal{G}$ which is the closest in direction to the straight line $\overline{OX}$. As pointed out by Liebeherr et al. in [20], on a Delaunay triangulation, compass routing is a proper navigation. The associated directed navigation with direction $e_1$ is

$$\mathcal{A}_{e_1}(X) = \arg\max\left\{\left\langle e_1, \frac{X-Y}{|X-Y|}\right\rangle : (X,Y) \in E\right\}.$$

*Radial navigation.* For $X, Y \in N^O$, $X \neq O$, $|Y| < |X|$, the radial navigation is defined by

$$\mathcal{A}(X) = |Y| \quad \text{if } N(B(X, |X-Y|) \cap B(O, |X|)) = \varnothing.$$

$\mathcal{A}(X)$ is the closest point to $X$ which is closer to the origin. Radial navigation has an a.s. positive progress and $\mathcal{A}(X)$ is a.s. uniquely defined. The corresponding navigation tree is the radial spanning tree and it is analyzed in [3]. The directed navigation associated with radial navigation is as follows: if $X, Y \in N$ and $\langle Y - X, e_1 \rangle > 0$, then

$$\mathcal{A}_{e_1}(X) = Y \quad \text{if } N(B(X, |X-Y|) \cap \mathcal{H}_{e_1}(X)) = \varnothing.$$

1.6. *Overview and organization of the paper.* In this paper, we illustrate our results on the small world navigation on $N^O$. In this paragraph, we state the main results that our analysis implies for this navigation.

*Local weak convergence of the navigation tree.* In Section 2, for the small world navigation, we prove that the navigation tree converges to a directed navigation forest for the local weak convergence on graphs as defined by Aldous and Steele [1]. More precisely, for a graph $\mathcal{G} = (N, E)$, we define $S_Y \circ \mathcal{G} = (S_Y N, E)$ as the graph obtained by translating all vertices $N$ by $Y$ and keeping the same edges. For the small world graph, let $\mathcal{T}_{e_1}(N)$ [resp. $\mathcal{T}_O(N)$] denote the directed navigation forest with direction $e_1$ (resp. navigation tree) built on the point set $N$.

PROPOSITION 1.3.  *Assume $\beta > d$ in the small-world navigation. If $|X|$ tends to $+\infty$ and $X/|X|$ to $e_1 \in S^{d-1}$, then $S_{-X} \circ \mathcal{T}_O(N^{O,X})$ converges to $\mathcal{T}_{-e_1}(N^O)$ for the local weak convergence.*

*Path average.* Let $g$ be a measurable function from $\mathbb{R}^d \times \mathbb{R}^d$ to $\mathbb{R}$ and $G(X) = \sum_{k=0}^{H(X)-1} g(X_k, X_{k+1})$, $G(O) = 0$. In Section 3, we will state various convergence results for $G(X)$ as $|X|$ tends to $\infty$ for a regenerative navigation.

*Proving that a navigation is regenerative.* Section 4 is the main contribution of this paper, wherein we prove that the small world navigation is regenerative. The method relies on the geometric properties of the navigation and tail bounds in the GI/GI/$\infty$ queue. It is based on a coupling which is related to the pseudo-regenerative times in Markov chains; see Athreya and Ney [2]. We prove that the small world navigation on $N^O$ has good regenerative properties for $\beta \leq d$ and $\beta > d+2$. We cannot extend this method to the case $d < \beta \leq d+2$. The method can, however, be extended to other navigation algorithms; see [7].

THEOREM 1.4. *For the small world navigation on $N^O$, if $\beta > d+1$ or $d-2 < \beta \leq d$, then $\mathcal{A}$ is regenerative. Moreover:*

- *if $\beta > d+2$, then there exists $\mu'$ such that a.s.*

$$\lim_{|X| \to \infty} \frac{H(X)}{|X|} = \frac{1}{\mu'};$$

- *if $\beta = d$, then there exists $\tilde{\mu}$ such that a.s.*

$$\lim_{|X| \to \infty} \frac{H(X)}{\ln |X|} = \frac{1}{\tilde{\mu}};$$

- *if $d-2 < \beta < d$, then a.s.*

$$\lim_{|X| \to \infty} \frac{H(X)}{\ln \ln |X|} = -\frac{1}{\ln(1-(d-\beta)/2)}.$$

The constant $\mu'$ is not computed explicitly; an expression will be given in the forthcoming Theorem 4.1. $\mu'$ is the asymptotic mean directed progress. We will exhibit an analytical expression for the constant $\tilde{\mu}$.

*Path deviation and tree topology.* In Section 5, we examine the path from $X$ to $O$ in the navigation tree. For regenerative navigation algorithms, we establish an upper bound on the maximal deviation of this path with respect to the straight line $\overline{OX}$:

(3) $$\Delta(X) = \max_{0 \leq k \leq H(X)} |X_k - \overline{X}_k|,$$

where $\overline{X}_k = \langle X_k, X/|X| \rangle X/|X|$ is the projection of $X_k$ on the straight line $\overline{OX}$.

An important feature of a tree is its set of *ends*. An end is a semi-infinite self-avoiding path in $\mathcal{T}_O$, starting from the origin: $(O = Y_0, Y_1, \ldots)$. The set of ends of a tree is the set of distinct ends (two semi-infinite paths are not distinct if they share an infinite subpath). A semi-infinite path $(O = Y_0, Y_1, \ldots)$ has an asymptotic direction if $Y_n/|Y_n|$ has a limit in the unit sphere $S^{d-1}$. Following Howard and Newman [15], some properties of the ends of $\mathcal{T}_O$ will follow from tail bounds on $\Delta(X)$. For $X \in N$, let $\Pi_{\text{out}}(X)$ be the set of offsprings of $X$ in $\mathcal{T}_O$, namely the set of points $Y \in N$ such that $X \in \pi(Y)$.

DEFINITION 1.5 (Howard and Newman [15]). Given $f \in \ell_0(\mathbb{R}_+)$, a tree is said to be $f$-straight at the origin, if, for all but finitely many vertices,

$$\Pi_{\text{out}}(X) \subset C(X, f(|X|)),$$

where, for all $X \in \mathbb{R}^d$ and $\epsilon \in \mathbb{R}^+$, $C(X,\epsilon) = \{Y \in \mathbb{R}^d : \theta(X,Y) \le \epsilon\}$ and $\theta(X,Y)$ is the angle (in $[0,\pi]$) between $X$ and $Y$.

[Recall that $\ell_0(\mathbb{R}_+)$ is the set of functions tending to 0 at $+\infty$.] $f$-straight trees have a simple topology described by Proposition 2.8 of [15] and restated in the following proposition.

PROPOSITION 1.6 (Howard and Newman [15]). *Let $\mathcal{T}$ be an $f$-straight spanning tree on a PPP. The following set of properties holds almost surely:*

- *every semi-infinite path has an asymptotic direction;*
- *for every $u \in S^{d-1}$, there exists at least one semi-infinite path with asymptotic direction $u$;*
- *the set of $u$'s of $S^{d-1}$ such that there is more than one semi-infinite path with asymptotic direction $u$ is dense in $S^{d-1}$.*

$f$-straightness will be related to $\Delta(X)$. On the small world navigation, we obtain the following proposition.

THEOREM 1.7. *For the small world navigation on $N^O$, there exists $C \ge 1$ such that if $\beta > (C+1)d + C$, then, for all $C(d+1)/(\beta - d) < \gamma < 1$, there exist $\eta > 0$ and $C_1 > 0$ such that*

$$\mathbb{P}(\Delta(X) \ge |X|^\gamma) \le C_1 |X|^{-d-\eta}$$

*and $\mathcal{T}_O$ is $f$-straight with $f(x) = |x|^{\gamma-1}$.*

*Shape of the navigation tree.* Finally, in Section 6, we will state a shape theorem for regenerative navigation algorithms. We define

$$\mathcal{T}_O(k) = \{X \in N : \mathcal{A}^k(X) = O\}.$$

On the small-world graph, we will obtain the following proposition.

THEOREM 1.8. *For the small-world navigation on $N^O$, let $\mu'$ and $\tilde\mu$ be as in Theorem 1.4. There exists $C \ge 1$ and a.s. for all $\epsilon > 0$ there exists $K \in \mathbb{N}$ such that the following all hold.*

- *If $\beta > (C+1)d + 2C$ and $k \ge K$, then*

$$N \cap B(O, (1-\epsilon)k\mu) \subset \mathcal{T}_O(k) \subset B(O, (1+\epsilon)k\mu).$$

*Moreover, a.s. and in $L^1$, $\frac{|\mathcal{T}_O(k)|}{\pi_d k^d} \to \mu'^d$.*

- If $\beta = d$ and $k \geq K$, then
$$N \cap B(O, e^{(1-\epsilon)k\tilde{\mu}}) \subset \mathcal{T}_O(k) \subset B(O, e^{(1+\epsilon)k\tilde{\mu}}).$$
Moreover, a.s. and in $L^1$, $\frac{\ln |\mathcal{T}_O(k)|}{k} \to d\tilde{\mu}$.
- For $d-2 < \beta < d$, let $\alpha = 1 - (d-\beta)/2$. If $k \geq K$, then
$$N \cap B(O, \exp(\alpha^{(1-\epsilon)k})) \subset \mathcal{T}_O(k) \subset B(O, \exp(\alpha^{(1+\epsilon)k})).$$
Moreover, a.s. and in $L^1$, $\frac{\ln \ln |\mathcal{T}_O(k)|}{k} \to \ln \alpha$.

## 2. Local convergence of navigation to directed navigation.

2.1. *Local weak convergence and proof of Proposition* 1.3. We now introduce *stable functionals* (see Lee [19] or Penrose and Yukich [24, 25]). For a Borel set $B$, $\mathcal{F}_B^N$ denotes the smallest $\sigma$-algebra such that the point set $N \cap B$ is measurable.

DEFINITION 2.1. Let $\mathcal{G} = (N, E)$ be a graph on $N$ and $F(X, \mathcal{G})$ be a measurable function valued in $\mathbb{R}$. $F$ is *stable* on $\mathcal{G}$ if, for all $X \in \mathbb{R}^d$, there exists a random variable $R > 0$ such that $F(X, \mathcal{G})$ is $\mathcal{F}_{B(X,R)}^N$-measurable and the distribution of $R$ does not depend on $X$.

Recall that $\mathcal{T}_{e_1}$ denotes the directed navigation forest associated with $\mathcal{A}_{e_1}$ on $N$ and $\mathcal{T}_O$ the navigation tree associated with $\mathcal{A}$ on $N^O$.

THEOREM 2.2. *(For the small-world navigation on a PPP and $\beta > d$.)* Let $F$ be a stable functional on $\mathcal{T}_{-e_1}$. As $x$ tends to $+\infty$, the distribution of $F(xe_1, \mathcal{T}_O)$ converges in total variation toward the distribution of $F(O, \mathcal{T}_{-e_1})$.

Proposition 1.3 is a corollary of Theorem 2.2. Using coupling, we start with an upper bound for the total variation distance between the distribution of $\mathcal{A}(X)$ (built on $N^{O,X}$) and $\mathcal{A}_{e_1}(X)$ (built on $N^X$). The proof of the next lemma is postponed to the Appendix.

LEMMA 2.3. *(For the small world navigation on a PPP and $\beta > d$.)* Let $X \in \mathbb{R}^d \setminus \{O\}$ and $e_1 \in S^{d-1}$ with $\cos \theta = \langle X/|X|, e_1 \rangle$. There exists a function $\varepsilon \in \ell_0(\mathbb{R}_+)$ and a coupling of $\mathcal{A}(X)$ and $\mathcal{A}_{-e_1}(X)$ such that

(4) $$\mathbb{P}(\mathcal{A}(X) \neq \mathcal{A}_{-e_1}(X)) \leq \varepsilon(|X|) + \varepsilon(1/|\theta|).$$

PROOF OF THEOREM 2.2. We set $X = xe_1$, $x > 0$ and build $\mathcal{T}_O$ on $N^{O,X}$ and $\mathcal{T}_{-e_1}$ on $N^X$. For all $t > 0$, we define the event $J_t(X) = \{\mathcal{T}_O \cap B(X, t) =$ $\mathcal{T}_{-e_1} \cap B(X,t)\}$. $F$ is a stable functional on $\mathcal{T}_{-e_1}$ with associated radius $R = R(X)$ and we have

$$\mathbb{P}(F(X,\mathcal{T}_O) \neq F(X,\mathcal{T}_{-e_1}))$$
$$\leq \mathbb{P}(J_{R(X)}(X)^c)$$
$$\leq \mathbb{P}(R > t) + P(J_t(X)^c)$$
$$\leq \mathbb{P}(R > t) + \mathbb{P}\left(\bigcup_{Y \in N \cap B(X,t)} \mathcal{A}(Y) \neq \mathcal{A}_{-e_1}(Y)\right)$$
$$\leq \mathbb{P}(R > t) + \mathbb{P}(N(B(X,t)) \geq n) + n\varepsilon\left((x-t) + \frac{x-t}{t}\right),$$

where we have used Lemma 2.3. It follows easily that $\lim_{x \to \infty} \mathbb{P}(F(xe_1,\mathcal{T}_O) \neq F(xe_1,\mathcal{T}_{-e_1})) = 0$. To complete the proof, note that $\mathcal{T}_{-e_1}$ is stationary: $F(xe_1,\mathcal{T}_{-e_1})$ and $F(O,\mathcal{T}_{-e_1})$ have the same distribution. $\square$

2.2. *Progress distribution in the small world.* We consider the small-world navigation $\mathcal{A}$ and the directed small world navigation $\mathcal{A}_{e_1}$. This directed navigation is defined if and only if $\beta > d$. Let $F$ denote the distribution function of the directed progress $P_{e_1}(X) = \langle \mathcal{A}_{e_1}(X) - X, e_1 \rangle$ (which does not depend on $e_1$ and $X$) and $F_X$ the distribution function of the progress at $X$ in the small world $P(X) = |X| - |\mathcal{A}(X)|$.

LEMMA 2.4. *For the small world navigation on a PPP and $d \geq 2$, the following properties hold:*

1. *if $\beta > d$, then, as $t$ goes to infinity,*

$$\overline{F}(t) \sim t^{d-\beta} \frac{2c\omega_{d-2}}{\beta - d}\left(1 - e^{-\int_{\mathcal{H}_{e_1}(O)} f(y)dy}\right)^{-1} \int_0^{\pi/2} \cos^{\beta-d}\theta\, d\theta;$$

2. *if $\beta > d$, then, for all function $\varepsilon \in \ell_0(\mathbb{R}_+)$,*

$$\lim_{|X| \to +\infty} \sup_{t \leq |X|\varepsilon(|X|)} t^{\beta-d}|F_X(t) - F(t)| = 0;$$

3. *if $d - 2 < \beta < d$, then the distribution of $|\mathcal{A}(X)|/|X|^{1-(d-\beta)/2}$ converges weakly and*

$$\sup_{X:|X| \geq 1} \mathbb{E}\left[\left|\ln \frac{|\mathcal{A}(X)|}{|X|^{1-(d-\beta)/2}}\right| \Big| \mathcal{A}(X) \neq 0\right] < \infty;$$

4. *if $\beta = d$, and $\tilde{F}_X$ is the distribution of $\tilde{P}(X) = -\ln(1 - P(X)/|X|) \in [0,+\infty]$, then $\tilde{F}_X$ converges weakly to a cumulative distribution function $\tilde{F}$ with $\int \tilde{F}(s)\,ds = \tilde{\mu} \in (0,+\infty)$ and, moreover, $(\tilde{P}_X \mathbb{1}(\tilde{P}_X < \infty))_{X \in \mathbb{R}^d}$ is uniformly integrable.*

The distribution $\tilde{F}$ in statement 4 is given by (38) and the weak limit of $|\mathcal{A}(X)|/|X|^{1-(d-\beta)/2}$ has a distribution obtained in (37). For $d \geq 3$ and $0 < \beta < d-2$, similar convergence results hold. To avoid longer computations, we will not state them. The proof is postponed to the Appendix.

**3. Path average for memoryless and regenerative navigation.** In this section, under various assumptions, we derive the asymptotic behavior of $H(X)$, the generation of $X$ in the navigation tree $\mathcal{T}_O$ when $\mathcal{A}$ is a regenerative navigation.

3.1. *Finite mean progress.*

PROPOSITION 3.1. *Let $\mathcal{A}$ be a memoryless navigation with nonnegative progress. Let $F_X(t) = \mathbb{P}(P(X) \geq t)$ and assume that $F_X$ converges weakly to $F$ as $|X|$ tends to $\infty$ and that $(F_X)_{X \in \mathbb{R}^d}$ is uniformly integrable. Then, a.s.*

$$\lim_{|X| \to +\infty} \frac{H(X)}{|X|} = \frac{1}{\mu},$$

*where $\mu = \int_0^\infty rF(dr) < \infty$.*

Before proving this proposition, we state a simple lemma.

LEMMA 3.2. *Let $\mathcal{A}$ be a navigation with nonnegative progress on a PPP. Let $x_0 \geq 0$, $\tau(X) = \inf\{k \geq 0 : |\mathcal{A}^k(X)| \leq x_0\}$ and $s(\cdot) \in \ell_\infty(\mathbb{R}^d)$. If a.s. (resp. in $L^p$) $\tau(X)/s(X)$ converges to $Z$, then a.s. (resp. in $L^p$) $H(X)/s(X)$ converges to $Z$.*

PROOF. The progress is a.s. positive: $\mathcal{A}(X) \in B(O, |X|)$. It follows that $\tau(X) \leq H(X) \leq \tau(X) + \sup_{Y \in B(O,x_0) \cap N} H(Y) \leq \tau(X) + N(B(O, x_0))$. We conclude by noticing that for $C > 0$, $\mathbb{E} \exp(CN(B(O, x_0))) < \infty$. □

PROOF OF PROPOSITION 3.1. Assume, first, that $\mu > 0$. Letting $0 < \eta < \mu/2$, we may find $x_0 \leq x_1$ and a function $h$ such that if $|X| \geq x_1$, then

$$\mathbb{1}(t \leq x_0)(\overline{F}(t) - h(t)) \leq \overline{F}_X(t) \leq \overline{F}(t) + h(t),$$

where $\int h(t)\,dt \leq \eta$, $h(t) \leq \overline{F}(t)$ and $\int_0^{x_0} \overline{F}(t) - h(t)\,dt \geq \mu - 2\eta$. Let $\tau(X) = \inf\{n : |X_n| \leq x_1\}$ and $(U_n), n \in \mathbb{N}$ [resp. $(V_n), n \in \mathbb{N}$] be an i.i.d. sequence of variables with tail distribution $1 \wedge (\overline{F} + h)$ (resp. $\overline{F} - h$). We define $Y_n = |X| - \sum_{k=0}^{n-1} U_k$, $Z_n = |X| - \sum_{k=0}^{n-1} V_k \mathbb{1}(V_k \leq x_0)$, $\tau_+(X) = \inf\{n : |Y_n| \leq x_1\}$, $\tau_-(X) = \inf\{n : |Z_n| \leq x_1\}$ and obtain

$$\mathbb{1}(\tau(X) > n)Z_n \leq_{st} \mathbb{1}(\tau(X) > n)|X_n| \leq_{st} \mathbb{1}(\tau(X) > n)Y_n.$$

We deduce that

$$\tau_-(X) \leq_{st} \tau(X) \leq_{st} \tau_+(X). \tag{5}$$

We have $\mathbb{E}U_n \leq \mu + \eta$ and $\mathbb{E}V_n \mathbb{1}(V_n \leq x_0) \geq \mu - 2\eta$. By the renewal theorem, a.s.

$$\liminf_X \frac{\tau_-(X)}{|X|} \geq \frac{1}{\mu + \eta} \quad \text{and} \quad \limsup_X \frac{\tau_+(X)}{|X|} \leq \frac{1}{\mu - 2\eta}. \tag{6}$$

By (5) and (6) we obtain a.s. $\liminf_X \tau(X)/|X| \geq 1/(\mu+\eta)$ and $\limsup_X \tau(X)/|X| \leq 1/(\mu-2\eta)$. Then, by Lemma 3.2, $H(X)/|X|$ tends a.s. to $1/\mu$. If $\mu = 0$, consider only $\tau_-(X)$. $\square$

REMARK 3.3. Using results on renewal processes, the case $\overline{F}(t) \sim_{t\to\infty} ct^{-\alpha}$, $\alpha \in (0,1)$, $c > 0$ is treated in [7] and the scaling obtained for $H(X)$ is $|X|^\alpha$ for $\alpha \in (0,1)$ and $|X|/\ln|X|$ for $\alpha = 1$.

3.2. *Scaled progress*. In this paragraph, we discuss cases when $(|X| - P(X))|X|^{-\alpha}$ converges for some $0 < \alpha \leq 1$.

3.2.1. *Scale-free progress*. (For a definition of scale-free navigation, see Franceschetti and Meester in [10].) Let $\tilde{P}(X) = -\ln(1 - P(X)/|X|) \in \mathbb{R}_+ \cup \{+\infty\}$ and $\tilde{F}_X(t) = \mathbb{P}(\tilde{P}(X) \leq t)$. Note that $\mathbb{P}(\tilde{P}(X) = \infty)$ may be positive.

PROPOSITION 3.4. *Let $\mathcal{A}$ be a memoryless navigation with nonnegative progress. If $\tilde{F}_X$ converges weakly to $\tilde{F}$ as $|X|$ tends to infinity and $(\tilde{P}_X \mathbb{1}(\tilde{P}_X < \infty))_{X \in \mathbb{R}^d}$ is uniformly integrable, then a.s.*

$$\lim_{|X| \to +\infty} \frac{H(X)}{\ln|X|} = \frac{1}{\tilde{\mu}},$$

*where $\tilde{\mu} = \int s\tilde{F}(ds) < \infty$.*

PROOF. Defining, for $0 \leq i < H(X) - 1$, $\tilde{P}_i = -\ln(1 - P(X_i)/|X_i|)$, we have $|X_k| = |X| \prod_{i=0}^{k-1}(1 - P(X_i)/|X_i|)$ and $\ln|X_k| = \ln|X| - \sum_{i=0}^{k-1} \tilde{P}_i$. The corresponding path in $\mathbb{R} \cup \{-\infty\}$ is $\tilde{\pi}(X) = \{\ln|X|, \ln|X| - \tilde{P}_0, \ldots, -\infty\}$. Let $\tau(X) = \sup\{n : \ln|X_n| < 0\}$. From Lemma 3.2, a.s. $\tau(X)$ and $H(X)$ are equivalent as $|X|$ tends to infinity (provided that they tend to infinity). We may apply Proposition 3.1 to the path $\{\ln|X|, \ldots, \ln|X_{\tau(X)}|\}$. $\square$

3.2.2. *Sublinear scaling.* We study the case when $Q(X) = |\mathcal{A}(X)||X|^{-\alpha} = (|X| - P(X))|X|^{-\alpha}$ has a nondegenerate limit for some $0 < \alpha < 1$.

PROPOSITION 3.5. *Let $\mathcal{A}$ be a memoryless navigation with nonnegative progress. Assuming that $\sup_{X:|X|\geq 1} \mathbb{E}[|\ln Q(X)||Q(X) \neq 0] < +\infty$, then a.s.*

$$\lim_{|X|\to+\infty} \frac{H(X)}{\ln\ln|X|} = -\frac{1}{\ln\alpha}.$$

PROOF. For $1 \leq k \leq H(X)$, let $Q_k = |X_k||X_{k-1}|^{-\alpha}$. If $k < H(X)$, then $\ln|X_k| = \alpha^k \ln|X| + \sum_{i=1}^{k} \alpha^{k-i} \ln Q_i$, hence

$$\ln|X_k| = \alpha^k \ln|X| + R_k, \tag{7}$$

with $|R_k| \leq Z_k = \sum_{i=1}^{k} \alpha^{k-i} |\ln Q_i|$. With the convention that $Z_k = 0$ for $k \geq H(X)$, $(Z_k, X_k), k \in \mathbb{N}$ is a Markov chain and

$$Z_{k+1} = \alpha Z_k + |\ln Q_{k+1}|.$$

Let $0 < \beta < 1 - \alpha$. By assumption, there exists $C_1$ such that $\sup_{X \in \mathbb{R}^d} \mathbb{E}\mathbb{1}(Q(X) \neq 0)|\ln Q(X)| \leq C_1$ (with the convention "$0 \times \infty = 0$"). It follows that

$$\mathbb{E}(\mathbb{1}(H(X) > k+1)(Z_{k+1} - Z_k)|Z_k = z) \leq -(1-\alpha)z + C_1 \tag{8}$$
$$\leq -\beta z + C_1 \mathbb{1}(z \in C),$$

with $C = \{z \in \mathbb{R}_+ : z \leq C_1/(1 - \alpha - \beta)\}$. Equation (8) is a geometric drift condition on a Markov chain (see (V4), page 371 in Meyn and Tweedie [21]). Let $K = \inf\{k \geq 1 : Z_k \in C\}$. By Theorem 15.2.5 in [21], for some $s > 0$,

$$\sup_{z \in C} \mathbb{E}[e^{s(K \wedge H(X))}|Z_0 = z] < \infty. \tag{9}$$

Set $x_0 = \exp(1 + C_1/(1 - \alpha - \beta))$. By Lemma 3.2, it is sufficient to show that a.s.

$$\lim_{|X|\to+\infty} \frac{\tau(X)}{\ln\ln|X|} = -\frac{1}{\ln\alpha},$$

where $\tau(X) = \inf\{k \geq 0 : |X_k| \leq x_0\}$. We fix $\epsilon > 0$ and let $(X^n), n \in \mathbb{N}$ be a sequence in $\mathbb{R}^d$ such that $|X^n|$ tends to infinity. We define $K(n, \epsilon) = \lfloor -(1+\epsilon)(\ln\ln|X^n|)/(\ln\alpha)\rfloor$ and $K'(n, \epsilon) = H(X) \wedge \inf\{k \geq K(n, \epsilon) : Z_k \in C\}$. From the Borel–Cantelli lemma and (9), a.s. for $n$ large enough, $K'(n, \epsilon/3) \leq 2K(n, \epsilon/3) \leq K(n, \epsilon)$. Therefore, for $n$ large enough, from (7), we have

$$\ln|X_{K(n,\epsilon)}| \leq \ln|X_{K'(n,\epsilon/3)}| \leq (\ln|X^n|)^{-\epsilon/3} + C_1/(1-\alpha-\beta) \leq \ln x_0$$

and it follows that a.s.

$$\limsup_X \frac{\tau(X)}{\ln\ln|X|} \leq -\frac{1}{\ln\alpha}.$$

The same computation can be carried out with $K(n, -\epsilon)$ to get a lower bound. □

3.3. *Average along a path.* We have thus far taken interest only in $H(X)$. More generally, we may also consider $G(X) = \sum_{i=0}^{H(X)-1} g(X_i, X_{i+1})$, where $g$ is a measurable function on $\mathbb{R}^d \times \mathbb{R}^d$. The proof of the next lemma is omitted since it is identical to the proof of Proposition 3.1.

LEMMA 3.6. *Let $\mathcal{A}$ be a memoryless navigation with nonnegative progress. Assume that $H(X)$ tends almost surely to infinity, that $(g(X, \mathcal{A}(X)))_{X \in \mathbb{R}^d}$ converges weakly as $|X|$ tends to infinity and that $(g(X, \mathcal{A}(X)))_{X \in \mathbb{R}^d}$ is uniformly integrable. Then, a.s.*

$$\lim_{|X| \to \infty} \frac{G(X)}{H(X)} = \lim_{|X| \to +\infty} \mathbb{E} g(X, \mathcal{A}(X)).$$

3.4. *Path average for regenerative navigation.* The next lemma is elementary, but nevertheless useful. Let $\mathcal{A}$ be a regenerative navigation, $\theta(X)$ its regenerative time and $X_k = \mathcal{A}^k(X)$. We define $\theta_0 = 0$, $\theta_{k+1} = \theta(X_{\theta_k})$ and $H^\theta(X) = \inf\{k \geq 0 : \mathcal{A}^{\theta_k}(X) = O\} = \inf\{k \geq 0 : (\mathcal{A}^\theta)^k(X) = O\}$.

LEMMA 3.7. *Let $\mathcal{A}$ be a regenerative navigation with nonnegative progress and regenerative time $\theta$. Assume that there exists a function $s(\cdot) \in \ell_\infty(\mathbb{R}^d)$ such that $\lim_{|X| \to \infty} H^\theta(X)/s(X) = 1/\mu$, $\mu > 0$. If $\lim_{|X| \to \infty} \mathbb{E}\theta(X) = \overline{\theta} < \infty$, then a.s.*

$$\lim_{|X| \to +\infty} \frac{H(X)}{s(X)} = \frac{\overline{\theta}}{\mu}.$$

PROOF. Note that $\theta_{H^\theta(X)-1} < H(X) \leq \theta_{H^\theta(X)}$, hence $\frac{\theta_{H^\theta(X)-1}}{s(X)} < \frac{H(X)}{s(X)} \leq \frac{\theta_{H^\theta(X)}}{s(X)}$. Letting $\tilde{\mathcal{A}} = \mathcal{A}^\theta$, we can apply Lemma 3.6 to $g(X, \tilde{\mathcal{A}}(X)) = \theta(X)$ to obtain that $\theta_{H^\theta(X)}/H^\theta(X)$ converges almost surely to $\overline{\theta}$. □

**4. Proof of Theorem 1.4.** The method of proof of Theorem 1.4 can be extended to a broader context and can be applied to, for example, radial navigation; see [7].

4.1. *Directed navigation on a small world.* In this paragraph, we prove that the small-world directed navigation is regenerative. We consider the model introduced in Section 1.5 with $\beta > d$. The maximal progress navigation from $X \in N$ with direction $e_1 \in S^{d-1}$ is $\mathcal{A}_{e_1}(X) = \arg\max\{\langle Y, e_1 \rangle : U(X,Y) \leq f(|X-Y|), Y \in N\}$. As above, $F$ is the distribution function of the directed progress $P_{e_1}(X) = \langle \mathcal{A}_{e_1}(X) - X, e_1 \rangle$ (which does not depend on $e_1$ and $X$), $X_k = \mathcal{A}_{e_1}^k(O)$ and $P_{e_1,k} = P_{e_1}(X_k)$.

THEOREM 4.1. *For the small-world directed navigation on a PPP and $\beta > d$:*

- *if $\beta > d + 1$, then $\mathcal{A}_{e_1}$ is regenerative for a stopping time $\theta$;*
- *if $\beta > d + 2$, then $\mathbb{E}\theta < \infty$ and a.s.*

$$\lim_{k \to +\infty} \frac{\langle \mathcal{A}_{e_1}^k(O), e_1 \rangle}{k} = \mu' = \frac{\mathbb{E}\langle \mathcal{A}_{e_1}^\theta(O), e_1 \rangle}{\mathbb{E}\theta}.$$

The regenerative time $\theta$ is built with the help of coupling techniques, $\theta$ is a stopping time on an enlarged filtration $\overline{\mathcal{F}}_k = \sigma((X_0, \omega_0), \ldots, (X_k, \omega_k))$, where $(\omega_k)$ are independent PPP's of intensity 1, independent of $N$. The remainder of the paragraph is devoted to the proof of Theorem 4.1.

PROOF OF THEOREM 4.1.
*Step* 1: *Coupling.* Let $N_X = N \cap \mathcal{H}_{e_1}(X) - X$. $N_O$ is a PPP of intensity 1 in $\mathcal{H}_{e_1}(O)$. However, due to the dependency structure, $N_{X_k}$ is not a PPP of intensity 1 in $\mathcal{H}_{e_1}(O)$. The idea of the proof is to enlarge the filtration $\mathcal{F}_k$ to build a sequence $(Z_k)$, $Z_k \in \mathbb{R}_+$ such that $\langle X_k, e_1 \rangle \leq Z_k$ and $N_{Z_k e_1}$ is a PPP of intensity 1. It can then easily be checked that if $\langle X_k, e_1 \rangle = Z_k$, then $\theta = k$ is a regenerative time. □

LEMMA 4.2. *Given $(X_0, \ldots, X_n)$, $N_{X_n}$ is a PPP of intensity $\lambda_n(x) = \prod_{k=0}^{n-1}(1 - f(|x + X_n - X_k|))$ in $\mathcal{H}_{e_1}(O)$.*

PROOF. The proof is by induction. Set $X_0 = O$ and assume that $N_O$ is a PPP of intensity $\lambda_0(x)\,dx$. It is sufficient to prove that $N_{X_1}$ is a PPP of intensity $\lambda_1(x) = (1 - f(|x + X_1|))\lambda_0(x + X_1)$ in $\mathcal{H}_{e_1}(O)$.

The set of neighbors of $O$ in $\mathcal{H}_{e_1}(O)$ is denoted by $V_O$. It is a thinning of $N_O$ and $V_O$ is a nonhomogeneous PPP on $\mathcal{H}_{e_1}(O)$ with intensity $f(|x|)\lambda_0(x)\,dx$, conditioned on $\{V_O \neq \varnothing\}$. In the next computation, $\tilde{V}_O$ will denote a PPP on $\mathcal{H}_{e_1}(O)$ of intensity $f(|x|)\lambda_0(x)\,dx$. If $A$ is a Borel set in $\mathcal{H}_{e_1}(O)$, then $\mathbb{P}(V_O(A) = 0) = \mathbb{P}(\tilde{V}_O(A) = 0)\mathbb{P}(\tilde{V}_O(\mathcal{H}_{e_1}(O) \cap A^c) > 0)/\mathbb{P}(\tilde{V}_O(\mathcal{H}_{e_1}(O)) > 0)$. We then write

$$\mathbb{P}(N(A) = k | A \subset \mathcal{H}_{e_1}(X_1))$$
$$= \mathbb{P}(\tilde{V}_O(A) = 0 | N(A) = k)\mathbb{P}(N(A) = k)/\mathbb{P}(\tilde{V}_O(A) = 0)$$
$$= \left(1 - \frac{\int_A f(|x|)\lambda_0(x)\,dx}{|A|}\right)^k \frac{|A|^k}{k!} e^{-|A| + \int_A f(|x|)\lambda_0(x)dx}$$
$$= \left(|A| - \int_A f(|x|)\lambda_0(x)\,dx\right)^k \frac{1}{k!} e^{-|A| + \int_A f(|x|)dx}.$$

Therefore, $N \cap \mathcal{H}_{e_1}(X_1)$ is a PPP of intensity $(1 - f(|x|))\lambda_0(x)$. Hence, $N_{X_1}$ is a PPP of intensity $\lambda_1(x) = (1 - f(|x + X_1|))\lambda_0(x + X_1)$ in $\mathcal{H}_{e_1}(O)$. □

By Lemma 4.2, far from $X_1$, the distributions $N_{X_1}$ and $N_O$ are close. We formalize this idea with the next lemma.

LEMMA 4.3. *There exists a random variable $Y_1 = Y(X) \geq \langle X_1, e_1 \rangle$ such that for $t \in \mathbb{R}_+$,*

(10) $$N_{Y_1 e_1} | (X_1, Y_1) \stackrel{\mathcal{L}}{=} N_O,$$

(11) $$\mathbb{P}(Y_1 - \langle X, e_1 \rangle > t) \leq C_1 t^{d-\beta}.$$

$N_{Y_1 e_1} | (X_1, Y_1) \stackrel{\mathcal{L}}{=} N_O$ is a shortened notation for stating that the conditional law of $N_{Y_1 e_1}$, given the pair $(X_1, Y_1)$, is equal to the law of $N_O$.

PROOF. Since $X_0 = X$, $N \cap \mathcal{H}_{e_1}(X_1)$ is a PPP of intensity $(1 - f(|x - X|)) \, dx$. Let $\tilde{V}(X)$ be a PPP with intensity $f(|X - x|) \, dx$ and independent of $N$. Then, $(\tilde{V}(X) + N) \cap \mathcal{H}_{e_1}(X_1)$ is a PPP of intensity $1$ on $\mathcal{H}_{e_1}(X_1)$. $\tilde{V}(X)$ is a.s. a finite point set and letting $\rho(X) = \inf\{r > 0 : \tilde{V}(X) \subset B(X, r)\}$, we have

$$\mathbb{P}(\rho(X) \geq t) = 1 - \exp\left(-\int_{B(O,t)^c} f(x) \, dx\right) \leq C_1 t^{d-\beta}.$$

We define

$$Y_1 = Y(X) = \langle X, e_1 \rangle + \max(P_{e_1}(X), \rho(X))$$

and (11) holds. If $A$ is a Borel set in $\mathcal{H}_{e_1}(Y(X)e_1)$, then $(\tilde{V}(X) + N)(A) = N(A)$. Since $(\tilde{V}(X) + N) \cap \mathcal{H}_{e_1}(X_1)$ is a PPP of intensity $1$, we deduce (10). □

Now, let $\tilde{V}(X_k), k \in \mathbb{N}^*$ be an independent sequence of PPP's, independent of $N_{X_k}$, with intensity $f(|X_k - x|) \, dx$. We define $\rho_k = \inf\{r > 0 : \tilde{V}(X_k) \subset B(X_k, r)\}$,

$$Y_k = \max(\langle X_k, e_1 \rangle, \langle X_{k-1}, e_1 \rangle + \rho_{k-1}) \quad \text{and} \quad Z_k = \max_{1 \leq l \leq k} Y_l.$$

Let $\overline{\mathcal{F}}_k$ be the $\sigma$-algebra generated by $((X_1, Y_1), \ldots, (X_k, Y_k))$. For the sake of simplicity, we will simply write $\mathcal{F}_k$ for $\overline{\mathcal{F}}_k$. Using the same argument as in Lemma 4.3,

(12) $$N_{Z_k} | \mathcal{F}_k \stackrel{\mathcal{L}}{=} N_O.$$

Note also that, since $\rho_k$ is independent of $N$, we have

(13) $$\mathbb{P}(Y_k - \langle X_k, e_1 \rangle \geq t | \mathcal{F}_k) \leq \mathbb{P}(\rho_k \geq t) \leq C_1 t^{d-\beta}.$$

We endow the set of point processes with the natural partial order relation: $N \leq_{st} N'$ if, for all Borel sets $A$ and $t \in \mathbb{N}$, $\mathbb{P}(N(A) \geq t) \leq \mathbb{P}(N'(A) \geq t)$.

LEMMA 4.4.

(14) $$N_{X_k}|\mathcal{F}_k \leq_{st} N_O$$

and, for some $C_0 > 0$:

(i) $\mathbb{P}(P_{e_1,k} \geq 1|\mathcal{F}_k) \geq C_0$;
(ii) $\mathbb{P}(Y_k = \langle X_k, e_1\rangle|\mathcal{F}_k) \geq C_0$.

PROOF. Equation (14) is a direct consequence of the fact that $N_{X_k}$ is a nonhomogeneous PPP of intensity $\prod_{l=0}^{k-1}(1 - f(|x + X_k - X_l|)) \leq 1$ (by Lemma 4.2). Statement (ii) is a consequence of (i) since $\mathbb{P}(Y_k = \langle X_k, e_1\rangle|\mathcal{F}_k) \geq \mathbb{P}(\rho_k \leq 1)\mathbb{P}(P_{e_1,k} \geq 1|\mathcal{F}_k) \geq C_0$. Assertion (i) stems from the fact that the progress is a.s. positive. Indeed, let $t > 0$ and $\lambda_k(x) = \prod_{i=0}^{k-1}(1 - f(|x + X_k - X_i|))$. By Lemma 4.2, we have

$$\mathbb{P}(P_{e_1,k} \geq t|\mathcal{F}_k) = (1 - e^{-\int_{\mathcal{H}_{e_1}(te_1)} f(|x|)\lambda_k(x)dx})(1 - e^{-\int_{\mathcal{H}_{e_1}(O)} f(|x|)\lambda_k(x)dx})^{-1}$$

$$\geq \left(\int_{\mathcal{H}_{e_1}(te_1)} f(|x|)\lambda_k(x)\,dx\right)\left(\int_{\mathcal{H}_{e_1}(O)} f(|x|)\lambda_k(x)\,dx\right)^{-1}$$

$$\geq \left(\int_{\mathcal{H}_{e_1}(te_1)} f(|x|)\,dx\right)\left(\int_{\mathcal{H}_{e_1}(O)} f(|x|)\,dx\right)^{-1}.$$

[We have used the inequality $(1 - e^{-u})/(1 - e^{-U}) \geq u/U$ for all $0 < u \leq U$. The last inequality following easily from the fact that $\lambda_k(|x|)$ is a nondecreasing function of $|x|$ in $\mathcal{H}_{e_1}(O)$.] □

Note that (13) and Lemma 4.4(ii) imply that there exists a variable $\sigma$ such that

(15) $\quad Y_k - X_k|\mathcal{F}_k \leq_{st} \sigma, \qquad \mathbb{P}(\sigma = 0) > 0 \quad \text{and} \quad \mathbb{P}(\sigma \geq t) \leq C_1 t^{d-\beta}.$

*Step* 2: *Regenerative time.* We define $W_n = Z_n - \langle X_n, e_1\rangle \geq 0$, $W_0 = 0$. With the convention that inf over an empty set is $+\infty$ and letting $\theta_0 = 0$, $\theta_{n+1} = \inf\{k > \theta_n : W_k = 0\}$, $\theta = \theta_1$, we have

$$W_n \leq_{st} \left(\max_{2 \leq i \leq n-1}\left(\sigma_{i-1} - \sum_{k=i-1}^{n-1} \tau_k\right)\right)^+,$$

where $(\sigma_k)_{k \in \mathbb{N}}$ is a sequence of i.i.d. copies of $\sigma$ and $(\tau_k)_{k \in \mathbb{N}}$ is a sequence of i.i.d. copies of $\tau$ with $\mathbb{P}(\tau = 1) = C_0$ and $\mathbb{P}(\tau = 0) = 1 - C_0$, as in Lemma 4.4, assertion (i). Hence, $W_n$ is upper bounded by the largest residual service time in a GI/GI/$\infty$ queue (see Appendix A.3). By Lemma A.1, for $\beta > d + 1$, $\theta$ is a.s. finite and for $\beta > d + 2$, $\mathbb{E}\theta < \infty$.

*Step* 3: *Embedded memoryless directed navigation.* Assume $\beta > d + 2$. From Step 2, there exists an increasing sequence $(\theta_n), n \in \mathbb{N}$ with $\theta_0 = 0$, $\theta_1 = \theta$ and such that $(\theta_{n+1} - \theta_n)_{n \in \mathbb{N}}$ is i.i.d. and $\mathbb{E}\theta < \infty$. We define

$$P^\theta_{e_1,k} = \langle X_{\theta_{k+1}} - X_{\theta_k}, e_1 \rangle = \sum_{n=\theta_k}^{\theta_{k+1}-1} P_{e_1,n}.$$

The sequence $(P^\theta_{e_1,k}), k \in \mathbb{N}$ is i.i.d. It remains to check that $\mathbb{E}P^\theta_{e_1,0} < \infty$. Note that

(16) $\quad \mathbb{P}(P_{e_1,n} > t | \mathcal{F}_n) \leq \mathbb{1}(t > W_n)\overline{F}(t - W_n) + \mathbb{1}(t \leq W_n).$

By Lemma 2.4, as $t$ tends to infinity, we have $\overline{F}(t) \sim C_1 t^{d-\beta}$. Moreover, by Lemma A.1, $W_k \leq_{st} M$, where $M$ is the stationary solution of the GI/GI/$\infty$ queue $\mathbb{P}(M > t) \leq C_1 t^{1+d-\beta}$. It follows that

$$\mathbb{E}P^\theta_{e_1,0} \leq \mathbb{E}\theta \mathbb{E}P_{e_1,0} + \mathbb{E}\sum_{n=0}^{\theta-1} W_n < \infty.$$

($\mathbb{E}\sum_{n=0}^{\theta-1} W_n < \infty$ is due to the cycle formula or Kac's formula; see, e.g. Section 3.1 of Baccelli and Brémaud [4].) Then, Theorem 4.1 follows from the strong law of large numbers.

4.2. *Proof of Theorem* 1.4. The proof is an extension of the ideas presented in the previous paragraph. We will only point out the differences with the proof of Theorem 4.1. We set $X_k = \mathcal{A}^k(X)$ and $P_k = P(X_k) = |X_k| - |X_{k+1}|$.

4.2.1. *Proof of Theorem* 1.4: $\beta > d + 1$.
*Step* 1: *Coupling.*

LEMMA 4.5. *There exists a random variable $0 \leq Y(X) \leq |X| - P(X)$ such that, for all Borel sets $A$ with $A \subset B(O, Y(X))$, $t \in \mathbb{N}$,*

(17) $\quad \mathbb{P}(N(A) = t | (X_1, Y(X))) = \mathbb{P}(N(A) = t),$

(18) $\quad \mathbb{P}(|X| - Y(X) \geq t) \leq C_1 t^{d-\beta}.$

PROOF. $N^O \cap B(O, X_1)$ is a Poisson point process of intensity $(1 - f(|X - x|))\,dx$ under its Palm version at $O$. The proof uses the same coupling as Lemma 4.3. Let $\tilde{V}(X)$ be a PPP with intensity $f(|X - x|)\,dx$ and independent of $N$. Since $\tilde{V}(X) \cap B(O, |X| - P(X))$ and $N \cap B(O, |X| - P(X))$ are independent, $(\tilde{V}(X) + N) \cap B(O, |X| - P(X))$ is a PPP of intensity 1 on $B(O, |X| - P(X))$. $\tilde{V}(X)$ is a.s. a finite point set. Letting $\rho(X) = \inf\{r \geq 0 : \tilde{V}(X) \subset B(X, r)\}$, for some $C_1 > 0$ (not depending on $X$), we have

$$\mathbb{P}(\rho(X) \geq t) = 1 - \exp\left(-\int_{B(O,t)^c} f(x)\,dx\right) \leq C_1 t^{d-\beta}.$$

We then define $Y(X) = (|X| - \max(P(X), \rho(X)))^+$. □

Let $\tilde{V}(X_k)$ be a PPP with intensity $f(|X_k - x|)\,dx$, independent of $N$, and $\rho_k = \inf\{r \geq 0 : \tilde{V}(X_k) \subset B(X_k, r)\}$. We define $Y_0 = |X|$ and $Y_k = (\min(|X_{k-1}| - \rho_k, |X_k|))^+$ and let $\overline{\mathcal{F}}_k$ denote the $\sigma$-algebra generated by $(X_1, Y_1), \ldots, (X_k, Y_k)$. Let $Z_0 = |X|$ and

$$Z_k = \min(Z_{k-1}, |X_{k-1}| - \rho_k, |X_k|) = \min_{1 \leq l \leq k} Y_l.$$

For $A \subset B(O, Z_k)$, we have

(19) $$\mathbb{P}(N(A) = t | \mathcal{F}_k) = \mathbb{P}(N(A) = t).$$

The next lemma is the analog of Lemma 4.4. The proof is omitted.

LEMMA 4.6. *For all Borel sets $A \subset B(O, |X_k|)$, $t \in \mathbb{N}$,*

(20) $$\mathbb{P}(N(A) \geq t | \mathcal{F}_k) \leq \mathbb{P}(N(A) \geq t)$$

*and for some $C_0 > 0$:*

  (i) *if $|X_k| \geq 1$, then $\mathbb{P}(P_k \geq 1 | \mathcal{F}_k) \geq C_0$;*
  (ii) *if $|X_k| \geq 1$, then $\mathbb{P}(X_{k+1} = Y_{k+1} | \mathcal{F}_k) \geq C_0$.*

Since $\rho_k$ and $X_k$ are independent, $\mathbb{P}(|X_k| - Y_k \geq t | \mathcal{F}_k) \leq \mathbb{P}(\rho_k \geq t) \leq C_1 t^{d-\beta}$. Thus, Lemma 4.6(ii) implies that there exists a variable $\sigma$ such that, if $|X| \geq 1$, then

(21) $(Y_k - X_k) | \mathcal{F}_k \leq_{st} \sigma$, $\quad \mathbb{P}(\sigma = 0) > 0 \quad$ and $\quad \mathbb{P}(\sigma \geq t) \leq C_1 t^{d-\beta}$.

*Step* 2: *Regenerative time*. We define $W_n = |X_n| - Z_n \geq 0$, $W_0 = 0$ and for $n \geq H(X)$, $W_n = 0$. Then,

(22) $$W_n \leq_{st} \left(\max_{2 \leq i \leq n-1}\left(\sigma_{i-1} - \sum_{k=i-1}^{n-1} \tau_k\right)\right)^+,$$

where $(\sigma_k)_{k \in \mathbb{N}}$ is a sequence of i.i.d. copies of $\sigma$ given in (21) and $(\tau_k)_{k \in \mathbb{N}}$ is a sequence of i.i.d. copies of $\tau$ with $\tau = 1$ with probability $C_0$ and 0 otherwise, as in Lemma 4.6.

By (22), $W_n$ is upper bounded by the largest residual service time in a GI/GI/$\infty$ queue (see Appendix A.3). Let $\tilde{W}_n$ be the right-hand side of (22) and let $\theta = \inf\{k \geq 1 : \tilde{W}_k = 0\}$. By Lemma A.1 (in the Appendix), if $\beta > d+1$, then $\theta$ is a.s. finite and if $\beta > d+2$, then $\mathbb{E}\theta < \infty$. By (19), $\theta$ is a regenerative time for the small-world navigation [for $|X| < 1$, we set $\theta(X) = H(X)$].

*Step* 3: *Embedded memoryless navigation.* $\mathcal{A}^\theta$ is a memoryless navigation and we define $P^\theta(X) = |X| - |X_\theta| = \sum_{k=0}^{\theta-1} P_k$,

(23) $\qquad \mathbb{P}(P_k \geq t | \mathcal{F}_k) \leq \overline{F}_X(t - W_k)\mathbb{1}(t \geq W_k) + \mathbb{1}(t < W_k),$

where $W_k \leq_{st} M$, $\mathbb{P}(M \geq t) \leq C_1 t^{1-(\beta-d)}$ and $\overline{F}_X(t) \leq C_1 t^{d-\beta}$. If $(U_k)_{k\in\mathbb{N}}$ denotes an i.i.d. sequence of variables such that $\mathbb{P}(U_k \geq t) = 1 \wedge C_1 t^{\beta-d}$ with $U_k$ independent of $\mathcal{F}_k$, then we obtain from (23) that $P^\theta(X) \leq_{st} Q = \sum_{k=0}^{\theta-1}(U_k + W_k)$. Therefore, $\mathbb{E}Q = \mathbb{E}\theta \mathbb{E}U + \mathbb{E}\sum_0^{\theta-1} W_k < \infty$ and it also follows that $(P^\theta(X))_{X \in \mathbb{R}^d}$ is uniformly integrable.

The last step is to identify $\lim_{|X|\to\infty} \mathbb{E}P^\theta(X)$. For the directed navigation with direction $e_1$, the same regenerative time $\theta$ was introduced. Theorem 4.1 gives $\mathbb{E}P^\theta_{e_1}(O) = \mu'\mathbb{E}\theta$. Moreover, $P^\theta_{e_1}(X)$ is a stabilizing functional of the small world directed navigation tree $\mathcal{T}_{e_1}$ and the distribution of $P^\theta_{e_1}(O)$ does not depend on $e_1$. Hence, from Theorem 2.2, $P^\theta(X)$ converges weakly to $P^\theta_{e_1}(O)$. Since $(P^\theta(X))_{X\in\mathbb{R}^d}$ is uniformly integrable, we obtain

$$\lim_{|X|\to\infty} \mathbb{E}P^\theta(X) = \mathbb{E}P^\theta_{e_1}(O) = \mu'\mathbb{E}\theta.$$

Thus, we can apply Proposition 3.1 and Lemma 3.7 and deduce that $H(X)/|X|$ tends a.s. to $1/\mu'$.

4.3. *Proof of Theorem* 1.4: $\beta = d$. We define the scaled progress $\tilde{P}_k = -\ln(1 - P_k/|X_k|)$. We have $\ln|X_k| = \ln|X| - \sum_{i=0}^{k-1} \tilde{P}_k$. Up to this scaling, the proof is analogous to the case $\beta > d$. We need to be careful with the event $\{\tilde{P}(X) = \infty\} = \{P(X) = |X|\}$: in this section, we use the convention "$\ln\frac{0}{0} = +\infty$."

*Step* 1: *Coupling.* We define $Y(X) = \min(|\mathcal{A}(X)|, \sup\{r \geq 0 : B(O,r) \cap \tilde{V}(X) = \varnothing\})$, where $\tilde{V}(X)$ is a PPP intensity $f(|x - X|)$, independent of everything else. As in Lemma 4.5, we obtain the following.

LEMMA 4.7. *For all Borel sets $A$ with $A \subset B(O, Y(X))$, $t \in \mathbb{N}$,*

$$\mathbb{P}(N(A) = t | (X_1, Y(X))) = \mathbb{P}(N(A) = t).$$

*Moreover, for all $X \in \mathbb{R}^d$,*

$$\mathbb{P}\left(\ln\frac{|\mathcal{A}(X)|}{Y(X)} \geq s\right) \leq C_1 \exp(-2s).$$

We define the sequence $(Y_k)_{k\in\mathbb{N}^*}$ and $(Z_k)_{k\in\mathbb{N}^*}$ as previously: $Y_k = \min((|X_{k-1}| - \rho_{k-1})^+, |X_k|)$ and $Z_k = \min_{1 \leq l \leq k} Y_l$, where $\rho_k = |X_k| - \sup\{t :$

$B(O,t) \cap \tilde{V}(X_k) = \varnothing\}$ and $\tilde{V}(X_k)$ is a PPP with intensity $f(|x - X_k|)$ and independent of everything else. Equation (19) still holds, and the analog of Lemma 4.6 is as follows.

LEMMA 4.8. *For all Borel sets $A \subset B(O, |X_k|)$, $t \in \mathbb{N}$,*
$$\mathbb{P}(N(A) \geq t | \mathcal{F}_k) \leq \mathbb{P}(N(A) \geq t)$$
*and for some $C_0 > 0$:*

(i) *if $|X_k| \geq e$, then $\mathbb{P}(\tilde{P}_k \geq 1 | \mathcal{F}_k) \geq C_0$;*
(ii) *if $|X_k| \geq e$, then $\mathbb{P}(X_{k+1} = Y_{k+1} | \mathcal{F}_k) \geq C_0$.*

We deduce that then there exists a r.v. $\sigma$ such that if $|X_k| \geq e$, then

$$(24) \quad \ln \frac{|X_{k+1}|}{Y_{k+1}} | \mathcal{F}_k \leq_{st} \sigma, \qquad \mathbb{P}(\sigma = 0) > 0 \quad \text{and} \quad \mathbb{P}(\sigma \geq s) \leq C_1 e^{-2s}.$$

*Step* 2: *Regenerative time.* We define $H'(X) = \inf\{k \geq 0, |X_k| \leq e$ or $Y_k = 0\}$. By Lemma 3.2, if $H'(X)/\ln|X|$ converges then $H(X)/\ln|X|$ converges, to the same limit. Now, for $n \leq H'(X)$, set $W_n = \ln(|X_n|/Z_n) \geq 0$ and we have

$$W_{n+1} = \max\left(W_n - \ln \frac{|X_n|}{|X_{n+1}|}, \ln \frac{|X_{n+1}|}{Y_{n+1}}\right).$$

It then follows that

$$(25) \qquad W_n \leq_{st} \left(\max_{2 \leq i \leq n-1}\left(\sigma_{i-1} - \sum_{k=i-1}^{n-1} \tau_k\right)\right)^+,$$

where $(\sigma_k)_{k \in \mathbb{N}}$ is a sequence of i.i.d. copies of $\sigma$ given in (24) and $(\tau_k)_{k \in \mathbb{N}}$ is a sequence of i.i.d. copies of $\tau$ with $\mathbb{P}(\tau \geq 1) \geq C_0$. $W_n$ is upper bounded by the largest residual service time in a GI/GI/$\infty$ queue (see Appendix A.3). Let $\theta$ be the first positive time at which the queue appearing on the left-hand side of (25) is empty. By Lemma A.1 (in the Appendix), $\theta$ is a.s. finite and, for some $C > 0$,

$$(26) \qquad \mathbb{E} \exp(C\theta) < \infty.$$

*Step* 3: *Embedded memoryless navigation.* We define $\tilde{P}^\theta(X) = \sum_{k=0}^{\theta-1} \tilde{P}_k$. From (26), (23) and Lemma 2.4, we deduce that $(\mathbb{1}(\tilde{P}^\theta(X) < \infty)\tilde{P}^\theta(X))_{X \in \mathbb{R}^d}$ is uniformly integrable. We assume for now that $\tilde{P}^\theta(X)$ converges weakly as $|X|$ tends to infinity. Define $\tilde{\mu}' = \lim_{|X| \to \infty} \mathbb{E}\tilde{P}^\theta(X)/\mathbb{E}\theta$. From Proposition 3.4, we obtain

$$\mathbb{P}\text{-a.s.} \qquad \lim_{|X| \to \infty} \frac{H(X)}{\ln|X|} = \frac{1}{\tilde{\mu}'}.$$

It remains to prove that $\tilde{P}^\theta(X)$ converges weakly as $|X|$ tends to infinity (since $\beta = d$, we cannot apply Theorem 2.2). As already pointed out $N \cap B(O, |X_1|)$ is a PPP of intensity $\lambda_{X_0}(y) \, dy = (1 - f(|X_0 - y|)) \, dy$. The set of neighbors of $X_1$ in $B(O, |X_1|)$ is a PPP of intensity $\lambda_{X_0}(y) f(|X_1 - y|) \, dy$ with an extra atom at $O$ with probability $f(|X_1|)$, conditioned on not being empty, hence

$$\mathbb{P}(\tilde{P}_1(X) \geq s | \mathcal{F}_1) = \frac{A(s)}{A(0)},$$

with $A(s) = 1 - (1 - f(|X_1|)) \exp(-\int_{B(O,|X_1|e^{-s})} f(|X_1 - y|) \lambda_{X_0}(y) \, dy)$. With the change of variable $z = y/|X_1|$ and $e_i = X_i/|X_i|$, we obtain

$$A(s) = 1 - (1 - f(|X_1|))$$
$$\times \exp\left(-\int_{B(O,e^{-s})} |X_1|^d f(|X_1||e_1 - z|)(1 - f(|X_1||e_0 e^{\tilde{P}_0} - z|)) \, dz\right).$$

Since $|X_1| = |X|e^{-\tilde{P}_0}$, $\mathbb{P}(\tilde{P}_0(X) + \tilde{P}_1(X) \geq s)$ is equal to

$$1 - \mathbb{E}(1 - f(|X|e^{-\tilde{P}_0}))$$
$$\times \exp\left(-\int_{B(O,e^{-s-\tilde{P}_0})} e^{-d\tilde{P}_0} |X|^d f(e^{-\tilde{P}_0}|X||e_1 - z|)\right.$$
$$\left. \times (1 - f(|X||e_0 - ze^{-\tilde{P}_0}|)) \, dz\right),$$

divided by the same expression with $s = 0$. Letting $|X|$ tend to infinity, we deduce that $\tilde{P}_0(X) + \tilde{P}_1(X)$ converges weakly to $\tilde{Q}_0 + \tilde{Q}_1$, where $(\tilde{Q}_k)_{k \in \mathbb{N}}$ is an i.i.d. sequence of variables with common distribution function $\tilde{F}$. Similarly for all $n \in \mathbb{N}$, $\sum_{k=0}^{n-1} \tilde{P}_k(X)$ converges weakly to $\sum_{k=0}^{n-1} \tilde{Q}_k$. Since the sequence $(\sum_{k=0}^{\theta-1} \tilde{P}_k(X))$ is uniformly integrable, we deduce that $\lim_{|X| \to \infty} \mathbb{E} \sum_{k=0}^{\theta-1} \tilde{P}_k(X) = \mathbb{E}\theta \mathbb{E}\tilde{Q}_1$ and it follows that $\tilde{\mu} = \tilde{\mu}'$.

4.3.1. *Proof of Theorem* 1.4: $d - 2 < \beta < d$. The proof follows from Proposition 3.5 and the argument used in the case $\beta = d$. Let $\alpha = 1 - (d - \beta)/2$. We define, for $1 \leq k \leq H(X)$, $Q_k = |X_k|/|X_{k-1}|^\alpha$ and $Q_k = 0$, for $k > H(X)$.

Let $\rho_k = |X_k| - \sup\{r : B(O, r) \cap \tilde{V}(X_k) = \varnothing\}$, where $\tilde{V}(X_k)$ is a PPP with intensity $f(|x - X_k|)$ and independent of everything else. We define the sequences $(Y_k)_{k \in \mathbb{N}^*}$ and $(Z_k)_{k \in \mathbb{N}^*}$ as usual: $Y_k = \min((|X_{k-1}| - \rho_{k-1})^+, |X_k|)$ and $Z_k = \min_{1 \leq l \leq k} Y_l$. Letting $s > 0$, we have

$$\mathbb{P}\left(\ln \frac{|X_1|}{|Y_1|} \geq s \Big| Y_1 > 0\right) \leq \mathbb{P}\left(\ln \frac{|X_1|}{(|X| - \rho_0)^+} \geq s \Big| Y_1 > 0\right)$$

$$\leq \mathbb{P}\bigg(\ln Q_1 - \ln \frac{(|X|-\rho_0)^+}{|X|^\alpha} \geq s \bigg| Y_1 > 0\bigg)$$

$$\leq \mathbb{P}(|\ln Q_1| \geq s | X_1 \neq 0)$$

$$+ \mathbb{P}\bigg(\bigg|\ln \frac{|X|-\rho_0}{|X|^\alpha}\bigg| \geq s \bigg| \rho_0 < |X|\bigg)$$

$$\leq C_1 \exp(-C_0 s^2).$$

Let $H'(X) = \inf\{k \geq 0 : |X_k| \leq 1 \text{ or } Y_k = 0\}$ and $W_n = \ln(|X_n|/Z_n) \geq 0$. The remainder of the proof is similar to Section 4.3, with obvious changes.

## 5. Navigation tree topology.

5.1. *Maximal deviation, tree topology and $f$-straightness.* Let $\mathcal{A}$ be a navigation on $N^O$. Recall that $X_k = \mathcal{A}^k(X)$, $P_k = |X_k| - |X_{k+1}|$, $H(X) = \inf\{k : \mathcal{A}^k(X) = O\}$ and $\mathcal{F}_k = \sigma\{X_0, \ldots, X_k\}$. The path from $X$ to $O$ in the navigation tree $\mathcal{T}_O$ is denoted by $\pi(X) = \{X_0, X_1, \ldots, O\}$. The maximal deviation $\Delta(X)$ of this path is defined by (3). The next result relates $f$-straightness to $\Delta(X)$. It is an extension of Lemma 2.7 of [15].

PROPOSITION 5.1. *Let $\mathcal{A}$ be a navigation on $N^O$ and $\mathcal{T}_O$ its navigation tree. Let $\gamma \in (0,1)$ and $\eta > 0$. If $\mathbb{P}(\Delta(X) \geq |X|^\gamma) \leq C_1 |X|^{-d-\eta}$ and $\sup_{X \in \mathbb{R}^d} \mathbb{E}|X - \mathcal{A}(X)|^r < \infty$ for some $r > (d+1)/\gamma$, then $\mathcal{T}$ is $f$-straight at the origin for $f(x) = |x|^{\gamma-1}$.*

PROOF. If $N = \sum_n \delta_{T_n}$, then let $K$ be the number of points $T_n$ of $N$ such that $\Delta(T_n) \geq |T_n|^\gamma$. From the Slivnyak–Campbell formula,

$$\mathbb{E}K = \omega_{d-1} \int_0^\infty P(\Delta(x) \geq |x|^\gamma) x^{d-1} \, dx$$

$$\leq \omega_{d-1} \int_0^\infty x^{d-1} \min(1, C_1 x^{-d-\eta}) \, dx < \infty.$$

We define $B_{\gamma,x} = \{\exists X \in N : |X| \leq 2x \text{ and } |X - \mathcal{A}(X)| > x^\gamma\}$. Using the inequalities $\mathbb{P}(N(B(O,x)) \geq t) \leq \exp(-t \ln(t/(e\pi_d x^d)))$ and $\sup_{X \in \mathbb{R}^d} \mathbb{E}|X - \mathcal{A}(X)|^r \leq C_1$, we have

$$\mathbb{P}(B_{\gamma,x}) \leq \mathbb{P}(N(B(O,2x)) \geq e^2 \pi_d 2^d x^d) + e^2 \pi_d 2^d x^d \frac{\mathbb{E}|X - \mathcal{A}(X)|^r}{x^{r\gamma}}$$

$$\leq \exp(-e^2 \pi_d 2^d x^d) + e^2 C_1 \pi_d 2^d x^{d-r\gamma}$$

$$\leq C_1 x^{d-r\gamma}.$$

From the Borel–Cantelli lemma, it follows that there is some finite random $x_0$ such that for $X \in N \setminus B(O, x_0)$, $|X - \mathcal{A}(X)| \leq |X|^\gamma$. The rest of the proof uses the same argument as Lemma 2.7 of [15] (with $1 - \delta$ replaced by $\gamma$). □

5.2. *Memoryless isotropic navigation.* For $e_1, e_2 \in S^{d-1}$, we define $\mathbb{U}(e_1, e_2) = \{R \in \mathbb{U} : R(e_2) = e_1\}$, where $\mathbb{U}$ is the orthogonal group of $\mathbb{R}^d$, that is, $\mathbb{U}(e_1, e_2)$ is the set of isometries which maps $e_2$ to $e_1$.

DEFINITION 5.2. A navigation $\mathcal{A}$ on $N^O$ is *isotropic* if, for all $e_1, e_2$ in $S^{d-1}$, $x \geq 0$ and $R \in \mathbb{U}(e_1, e_2)$, $R\mathcal{A}(xe_2) \stackrel{\mathcal{L}}{=} \mathcal{A}(xe_1)$.

Let $X \neq O$ and $e_1, e_2$ in $S^{d-1}$, with $\langle e_i, X/|X|\rangle = 0$ for $i = 1, 2$. If $\mathcal{A}$ is isotropic, then $\langle \mathcal{A}(X), e_1\rangle \stackrel{\mathcal{L}}{=} \langle \mathcal{A}(X), e_2\rangle$. In particular, with $e_1 = -e_2$, we obtain $\mathbb{E}\langle \mathcal{A}(X), e_1\rangle = 0$. All examples in Section 1.5 are isotropic.

THEOREM 5.3. *Let $\mathcal{A}$ be a navigation on $N^O$. If $\mathcal{A}$ is isotropic and memoryless, with nonnegative progress, and:*

(i) $\sup_{X \in \mathbb{R}^d} \mathbb{E}|X - \mathcal{A}(X)|^r < \infty$, *with* $r > (d+1)/\gamma$;
(ii) *for* $|X| \geq x_0$, $\mathbb{P}(P(X) \geq C_0) \geq \epsilon$, *with* $x_0, C_0, \epsilon > 0$,

*then, for all $\eta \in (0, r - (d+1)/\gamma)$, there exists $C_1$ such that*

$$\mathbb{P}(\Delta(X) \geq |X|^\gamma) \leq C_1 |X|^{-d-\eta}$$

*and $\mathcal{T}_O$ is $f$-straight with $f(x) = |x|^{\gamma - 1}$.*

The straightness of $\mathcal{T}_O$ is a consequence of Proposition 5.1. As an immediate corollary, we have the following.

COROLLARY 5.4. *Under the assumption of Theorem 5.3, assume, further, that for all $r \geq 1$, $\sup_{X \in \mathbb{R}^d} \mathbb{E}|X - \mathcal{A}(X)|^r < \infty$. Then, for all $\eta > 0$ and $n \in \mathbb{N}$, there exists $C_1$ such that*

$$\mathbb{P}(\Delta(X) \geq |X|^{1/2+\eta}) \leq C_1 |X|^{-n}$$

*and $\mathcal{T}_O$ is $f$-straight with $f(x) = |x|^{-1/2+\eta}$.*

5.3. *Proof of Theorem 5.3.* Let $e_1, e_2 \in S^{d-1}$, with $\langle e_1, e_2\rangle = 0$, and assume that $X = |X|e_1$. We define

$$U_k = \langle X_k, e_1\rangle \quad \text{and} \quad V_k = \langle X_k, e_2\rangle.$$

Let $F = vect(e_1, e_2)$, $X_k^F$ be the orthogonal projection of $X_k$ on $F$ and $(\cos\theta_k, \sin\theta_k)$ be the coordinates of the projection of $X_k^F/|X_k^F|$ on the basis $(e_1, e_2)$. Let $R_k \in \mathbb{U}$ be such that $R_k X_k = |X_k|e_1$ and $e_2^k = R_k^{-1} e_2$. We define $p_k = \langle X_k^F - X_{k+1}^F, X_k/|X_k|\rangle$ and $q_k = \langle X_{k+1}^F - X_k^F, e_2^k\rangle$ (see Figure 1).

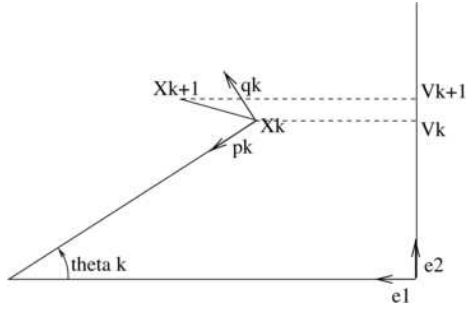

Fig. 1. $q_k$, $p_k$, $\theta_k$ and $V_k$.

We have $p_k \geq 0$ [since the navigation has nonnegative progress, $X_{k+1} \in B(O, |X_k|)$] and

(27)
$$V_{k+1} = V_k + q_k \cos\theta_k - p_k \sin\theta_k,$$
$$U_{k+1} = U_k - p_k \cos\theta_k - q_k \sin\theta_k,$$
$$\tan\theta_k = V_k/U_k.$$

If the navigation is isotropic and memoryless, then the distribution of $p_k$ and $q_k$ depends only on $|X_k|$, and $\mathbb{E}(q_k \cos\theta_k)^{2n+1} = 0$. We define

$$Q_k = q_k \cos\theta_k, \qquad S_k = \sup_{0 \leq n < k} |X_{n+1} - X_n| \quad \text{and} \quad M_k = \sup_{0 \leq \ell \leq k} \sum_{n=\ell}^{k-1} Q_n$$

(where, by convention, the sum over an empty set is equal to 0).

LEMMA 5.5. *If $0 \leq k \leq H(X)$, then $V_k \leq S_k + M_k$.*

PROOF. We prove this result by iteration. Assume that the inequality holds for $k-1$ and $\theta_{k-1} \in (0, \pi)$. Then, by (27), $V_k \leq V_{k-1} + Q_{k-1} \leq S_{k-1} + M_k$. Otherwise, $\theta_{k-1} \in (-\pi, 0]$, $V_{k-1} \leq 0$ and $V_k \leq S_k + V_{k-1} \leq S_k$. □

LEMMA 5.6. *Let $r' < r$. For all $t > 0$, there exists $C_t > 0$ such that $\mathbb{P}(\sup_{0 \leq \ell \leq k} \sum_{n=0}^{\ell-1} Q_n \geq k^\gamma t) \leq C_t k^{1-\gamma r'}$.*

PROOF. This lemma is a consequence of Theorem 3.1, equation (3.3) of Gut [14] (see also Theorem 2 in Baum and Katz [6]). This theorem is stated for a sum of independent variables, but it applies here as well. Indeed, note that:

1. $(|Q_k|)_{k \in \mathbb{N}}$ is stochastically bounded by an i.i.d. sequence $(Z_k)_{k \in \mathbb{N}}$ with $\mathbb{P}(Z_k \geq t | \mathcal{F}_k) \leq C_1 t^{-r}$, thus $\mathbb{E} Z^{r'} < \infty$;

2. $(Q_k)_{k \in \mathbb{N}}$ are nearly independent: if $n, m$, $k \ne l \in \mathbb{N}$, then $\mathbb{E}Q_k^{2n+1}Q_l^m = 0$ and $\mathbb{E}|Q_k^{2n}Q_l^{2m}| \le \mathbb{E}Z^{2n}\mathbb{E}Z^{2m}$.

Since the proofs of Gut, Baum and Katz rely only on the Markov inequality and truncation, their proofs apply here as well. □

LEMMA 5.7. *Letting $\mathcal{A}$ be as in Theorem 5.3, there exist positive constants $C_1$ and $C_0$ such that $\mathbb{P}(H(X) \ge |X|/\mu) \le C_1 \exp(-C_0|X|)$.*

PROOF. Let $\tau(X) = \inf\{k : |X_k| \le x_0\}$ and $(B_k), k \in \mathbb{N}$ be an i.i.d. sequence of Bernoulli variables with $\mathbb{P}(B_1 = 0) = 1 - \epsilon$ and $\mathbb{P}(B_1 = 1) = \epsilon$ [$\epsilon$, $c$ and $x_0$ were defined in assumption (ii) of Theorem 5.7]. Letting $\mu' < c\epsilon$, $\eta > 0$ and $\mu$ be such that $1/\mu = 1/\mu' + \eta$, we have

$$\mathbb{P}(H(X) \ge |X|/\mu) \le \mathbb{P}(N(B(O, x_0)) \ge \eta|X|) + \mathbb{P}(\tau(X) \ge |X|/\mu')$$

$$\le \mathbb{P}(N(B(O, x_0)) \ge \eta|X|) + \mathbb{P}\left(\sum_{k=0}^{\lfloor X/\mu' \rfloor - 1} P_k < |X| - x_0\right)$$

$$\le \mathbb{P}(N(B(O, x_0)) \ge \eta|X|) + \mathbb{P}\left(\sum_{k=0}^{\lfloor X/\mu' \rfloor - 1} B_k < |X|/C_0\right)$$

$$\le C_1 \exp(-C_0'|X|),$$

where we have used the inequality $\mathbb{P}(N(B(O, x_0)) \ge t) \le \exp(-t\ln(t/(e\pi_d x_0^d)))$ and Hoeffding's inequality: for $t < n\epsilon$, $\mathbb{P}(\sum_{k=0}^{n-1} B_k < t) \le 2\exp((t - n\epsilon)^2/(2n))$. □

We now conclude the proof of Theorem 5.3. Let $t > 0$ and $m_k = -\min_{0 \le \ell \le k} \times \sum_{n=\ell}^{k-1} Q_n$. By symmetry, Lemma 5.5 implies that $|V_k| \le S_k + M_k + m_k$, hence

$$\mathbb{P}\left(\max_{0 \le k \le H(X)} |V_k| \ge tx^\gamma, H(X) \le n\right)$$

$$\le 2\mathbb{P}\left(\max_{0 \le k \le n} M_k \ge tx^\gamma/3\right) + \mathbb{P}(S_n \ge tx^\gamma/3).$$

Note that $\max_{0 \le k \le n} M_k = \max_{0 \le k \le n} \sum_{\ell=0}^{k-1} Q_\ell - \min_{0 \le k \le n} \sum_{\ell=0}^{k-1} Q_\ell$ and so Lemma 5.6 gives $\mathbb{P}(\max_{0 \le k \le n} M_k \ge tx^\gamma/3) \le C_t n^{1-\gamma r'}$. We deduce that

$$\mathbb{P}\left(\max_{0 \le k \le H(X)} |V_k| \ge tx^\gamma, H(X) \le n\right) \le C_t n^{1-\gamma r'} + C_t nx^{-\gamma r'}.$$

From the isotropy of the navigation, we obtain

$$\mathbb{P}(\Delta(X) \ge x^\gamma)$$

$$\leq d\mathbb{P}\left(\max_{0\leq k\leq H(X)}|V_k|\geq x^\gamma/d\right)$$

$$\leq d\mathbb{P}\left(\max_{0\leq k\leq H(X)}|V_k|\geq x^\gamma/d, H(X)\leq x/\mu\right)+d\mathbb{P}(H(X)\geq x/\mu)$$

$$\leq C_1 x^{1-\gamma r'}+C_1\exp(-C_0 x).$$

If $r'<r$ is close enough to $r$, then $\gamma r'-1>d$ and this concludes the proof of Theorem 5.3.

5.4. *Isotropic regenerative navigation and proof of Theorem 1.7.* Let $\mathcal{A}$ be a regenerative navigation on $N^O$ with regenerative time $\theta$ so that $\mathcal{A}^\theta$ is a memoryless navigation. We define

$$L(X)=\sum_{k=0}^{\theta-1}|X_{k+1}-X_k|.$$

[Recall that $X_k=\mathcal{A}^k(X)$.] We have the following corollaries of Theorem 5.3.

COROLLARY 5.8. *Let $\mathcal{A}$ be a navigation on $N^O$ and $\gamma\in(1/2,1)$. If $\mathcal{A}$ is isotropic, regenerative and with nonnegative progress, and:*

(i) $\sup_{X\in\mathbb{R}^d}\mathbb{E}L(X)^r<\infty$, *with* $r>(d+1)/\gamma$;
(ii) *for* $|X|\geq x_0$, $\mathbb{P}(P(X)\geq c)\geq\epsilon$, *with* $x_0,c,\epsilon>0$;

*then the conclusions of Theorem 5.3 hold true.*

COROLLARY 5.9. *Under the assumption of Corollary 5.8, assume, further, that for all $r\geq 1$, $\sup_{X\in\mathbb{R}^d}\mathbb{E}L(X)^r<\infty$. The conclusions of Corollary 5.9 then hold true.*

PROOF OF COROLLARY 5.8. Let $\tilde{\mathcal{A}}(X)=\mathcal{A}^\theta(X)$, $\tilde{X}_k=\tilde{\mathcal{A}}^k(X)$, $\tilde{H}(X)$ be the length and $\tilde{\Delta}(X)$ be the deviation of the path $\{X,\tilde{X}_1,\ldots,O=\tilde{X}_{\tilde{H}(X)}\}$. Note that Theorem 5.3 applies to $\tilde{\mathcal{A}}$ and $\tilde{\Delta}(X)$. Therefore, Corollary 5.8 follows easily from Lemma 5.7 and the inequality, for all $0<r'<r$,

$$\mathbb{P}(\Delta(X)\geq 2|X|^\gamma)\leq\mathbb{P}(\tilde{\Delta}(X)\geq|X|^\gamma)$$

$$+\mathbb{P}\left(\max_{0\leq k\leq\lfloor|X|/\mu\rfloor}L(\tilde{X}_k)\geq|X|^\gamma\right)+\mathbb{P}(\tilde{H}(X)\geq|X|/\mu)$$

$$\leq C_1|X|^{1-\gamma r'}.\qquad\square$$

As an application, for the small-world navigation, we obtain Theorem 1.7. Indeed, for $\beta>d+2$, the small world navigation is isotropic regenerative with nonnegative progress. Moreover, we have $\mathbb{P}(|X-\mathcal{A}(X)|\geq t)\leq C_1 t^{d-\beta}$ and the tail of $\theta$ is bounded by a constant times $t^{2+d-\beta}$.

## 6. Shape of the navigation tree and proof of Theorem 1.8.

6.1. *Shape of memoryless navigation.* Let $\mathcal{A}$ be a navigation and let $\mathcal{T}_O(k) = \{X \in N : \mathcal{A}^k(X) = O\}$, the set of points at tree-distance less than $k$ from the origin.

THEOREM 6.1. *Let $\mathcal{A}$ be a memoryless navigation with nonnegative progress on $N^O$. Assume:*

(i) $\sup_{X \in \mathbb{R}^d} \mathbb{E} P(X)^r < \infty$ *for some* $r > d + 2$;
(ii) $F_X$ *converges weakly to $F$ with* $\mu = \int r F(dr) > 0$.

*Then, for all $\epsilon > 0$, there exists a.s. $K$ such that if $k \geq K$, then*

$$(28) \qquad N \cap B(O, (1-\epsilon)k\mu) \subset \mathcal{T}_O(k) \subset B(O, (1+\epsilon)k\mu).$$

*Moreover, a.s. and in $L^1$,*

$$(29) \qquad \frac{|\mathcal{T}_O(k)|}{\pi_d k^d} \to \mu^d.$$

In the literature, the constant $\mu$ is known as the *volume growth*. The intuition behind Theorem 6.1 is as follows: from Proposition 3.1, a point $k$ jumps away from the origin is asymptotically at Euclidean distance $D_k \sim k\mu$ from the origin. The ball of radius $D_k$ contains $\pi_d D_k^d$ points in $N$, asymptotically. In order to prove Theorem 6.1, we need an estimate of the tail of $H(X)$ around its mean.

PROPOSITION 6.2. *Under the assumption of Theorem 6.1, let $r' < r$. For all $\nu < \mu$, there exists a positive constant $C_1$ such that*

$$\text{if } |X| < n\nu - 1 \text{ and } n \geq 1, \qquad \text{then } \mathbb{P}(H(X) > n) \leq C_1 n(n\nu - |X|)^{-r'}.$$

*Similarly, for $\nu > \mu$,*

$$\text{if } |X| > 1 + n\nu \text{ and } n \geq 1, \qquad \text{then } \mathbb{P}(H(X) < n) \leq C_1 n(|X| - n\nu)^{-r'}.$$

In particular, if $\nu > \mu$, consider $n = \lfloor x/(2\nu) \rfloor$. We obtain

$$(30) \qquad P\left(H(X) > \frac{x}{\nu}\right) \leq C_1 |X|^{1-r'}$$

(and similarly for $\nu < \mu$).

Similarly, if assumption (i) is replaced by (i') $\sup_{X \in \mathbb{R}^d} \mathbb{E} \exp(sP(X)) < \infty$ for some $s > 0$, then we may prove an exponential tail bound for $H(X)$ (see [7]).

PROOF OF THEOREM 6.1. We define $G_k = |\mathcal{T}_O(k)| = \sum_{X \in N} \mathbb{1}(H(X) \leq k)$. $G_k$ is the size of the ball of center $O$ and radius $k$ for the graph-distance in $\mathcal{T}_O$. We start with the proof of (29). Letting $\epsilon \in (0,1)$, we write

$$|G_k - N^O(B(O, \mu k))| \leq \sum_{X \in N} \mathbb{1}(X \notin B(O, \mu k) \cap H(X) \leq k)$$
$$+ \sum_{X \in N} \mathbb{1}(X \in B(O, \mu k) \cap H(X) > k)$$
$$\leq \sum_{X \in N} \mathbb{1}(X \notin B(O, (1+\epsilon)\mu k) \cap H(X) \leq k)$$
$$+ N(B(O, (1+\epsilon)\mu k) \setminus B(O, (1-\epsilon)\mu k))$$
$$+ \sum_{X \in N} \mathbb{1}(X \in B(O, (1-\epsilon)\mu k) \cap H(X) > k)$$
$$\leq I_k + J_k + L_k.$$

From the Slyvniak–Campbell formula, using (30) for $\nu = \mu(1 + \epsilon/2)$, we obtain

$$\mathbb{E} I_k = \int_{\mathbb{R}^d \setminus B(O, (1+\epsilon)\mu k)} \mathbb{P}(H(X) \leq k)\, dX$$
$$\leq \omega_{d-1} \int_{(1+\epsilon)\mu k}^{\infty} C_1 x^{1-r'} x^{d-1}\, dx$$
$$\leq C_1 k^{d-r'+1}.$$

From the Borel–Cantelli lemma, we obtain that almost surely $I_k = 0$ for $k$ large enough. Similarly, letting $\nu = (1-\epsilon/2)\mu$, we obtain

$$\mathbb{E} L_k = \int_{B(O, (1-\epsilon)\mu k)} \mathbb{P}(H(X) \geq k)\, dX$$
$$\leq \omega_{d-1} \int_0^{(1-\epsilon)\mu k} C_1 k (k\nu - x)^{-r'} x^{d-1}\, dx$$
$$\leq \frac{\omega_{d-1} C_1 2^{r'}}{(k\epsilon\mu)^{r'-1}} \int_0^{(1-\epsilon)\mu k} x^{d-1}\, dx$$
$$\leq C_1 \epsilon^{1-r'} k^{d-r'+1}.$$

We deduce that almost surely $L_k = 0$ for $k$ large enough. Then, the ergodic properties of the PPP imply that

$$\frac{J_k}{k^d} = \frac{N(B(O, (1+\epsilon)\mu k) \setminus B(O, (1-\epsilon)\mu k))}{k^d}$$

converges almost surely and in $L^1$ toward $2d\pi_d(\mu\epsilon)^{d-1}$. We have thus proved that for all $\epsilon > 0$, almost surely

$$\limsup_k \frac{|G_k - N(B(O, \mu k))|}{k^d} \leq 2d\pi_d(\mu\epsilon)^{d-1}.$$

Hence, almost surely

$$\lim_k \frac{G_k}{k^d} = \lim_k \frac{N(B(O, \mu k))}{k^d} = \pi_d \mu^d.$$

The convergence in $L^1$ is a consequence of Scheffe's lemma. Equation (28) holds since a.s. for $k$ large enough, $I_k$ and $L_k$ are both equal to 0. $I_k$ is the cardinality of $\mathcal{T}_O(k) \cap B(O, (1+\epsilon)\mu k)^c$ and $L_k$ is the cardinality of $\mathcal{T}_O(k)^c \cap B(O, (1-\epsilon)\mu k)$.

6.2. *Proof of Theorem* 6.2. The following lemma is proved as Theorem 4 in Baum and Katz [6]. For details, see [7].

LEMMA 6.3. *Let* $(Y_k), k \in \mathbb{N}$ *be a sequence of i.i.d. real-valued random variables. We assume that* $\mathbb{E}Y_1 = 0$ *and* $\mathbb{E}|Y_1|^r < \infty$ *for some* $r > 1$. *Then, for all* $1 < r' < r$ *and* $t_0 \geq 0$, *there exists* $C_1$ *such that, for all* $t > t_0$ *and* $n \geq 1$,

$$\mathbb{P}\left(\left|\sum_{k=1}^n Y_k\right| \geq tn\right) \leq C_1 t^{-r'} n^{1-r'}.$$

As usual, for $1 \leq k \leq H(X)$, we define the progress $P_k(X) = |X_{k-1}| - |X_k|$ and for $k \geq H(X)$, $P_k(X) = 0$. We fix $r' < r'' < r$.

*Case* $\nu < \mu$. There exists $\nu' > 0$ such that $\nu' < \nu < \mu$ and $|X| < \nu' n - 1$. Since $(P(X))$ is uniformly integrable, there exists $x_0$ such that

(31) \qquad if $|X| \geq x_0$, \qquad then $\mathbb{E}P(X) \geq \nu$.

Letting $l < n$, we have

$$\mathbb{P}(H(X) > n)$$
$$= \mathbb{P}\left(\sum_{k=0}^{n-1} P_k < |X|\right)$$
$$\leq \mathbb{P}(N(B(O, x_0)) > l) + \mathbb{P}\left(\sum_{k=0}^{n-l-1} P_k < |X| - x_0\right)$$

(32) $\qquad \leq \exp\left(-l \ln \frac{l}{e\pi^d x_0^d}\right)$

$$+ \mathbb{P}\left(\sum_{k=0}^{n-l-1} P_k - \mathbb{E}(P_k|\mathcal{F}_k) < |X| - x_0 - \sum_{k=0}^{n-l-1} \mathbb{E}(P_k|\mathcal{F}_k)\right)$$

$$\leq \exp\left(-l \ln \frac{l}{e\pi^d x_0^d}\right)$$

$$+ \mathbb{P}\left(\sum_{k=0}^{n-l-1} P_k - \mathbb{E}(P_k|X_k) < |X| - (n-l)\nu\right),$$

where, in the last inequality, we have used (31). We define $Y_k = P_k - \mathbb{E}(P_k|X_k)$. By assumption (i),

$$\mathbb{E}Y_k = 0 \quad \text{and} \quad \mathbb{E}|Y_k|^r \leq C_1.$$

The sequence $(Y_k)_{k \in \mathbb{N}}$ is not independent, however:

1. $(|Y_k|)_{k \in \mathbb{N}}$ is stochastically dominated by an i.i.d. sequence $(Z_k)_{k \in \mathbb{N}}$ with $\mathbb{E}Z_k^{r''} < \infty$;
2. if $p, q \in \mathbb{N}$ and $k \neq l$, then $\mathbb{E}Y_k^{2p+1}Y_l^q = 0$ and $\mathbb{E}Y_k^{2p}Y_l^{2q} \leq \mathbb{E}Z_k^{2p}\mathbb{E}Z_l^{2q}$.

We can apply Lemma 6.3 to $(Y_k)$. It is stated for i.i.d. variables, but also holds in this case since the proof relies only on truncation and the Markov inequality. If $m \geq 1$ and $\nu - |X|/m > t_0$, $t_0 > 0$, then we obtain

$$\mathbb{P}\left(\sum_{k=0}^{m-1} P_k - \mathbb{E}(P_k|X_k) < |X| - m\nu\right) \leq C_1 m(m\nu - |X|)^{-r'}.$$

Hence, using this last inequality in (32) and setting $l = \lfloor(\nu'/\nu - 1)n\rfloor$, we obtain [since $(n-l)\nu \geq n\nu' > |X|$]

(33) $$\mathbb{P}(H(X) > n) \leq \exp(-C_0 n) + C_1 n(n\nu' - |X|)^{-r'}.$$

Then since $n \geq (n\nu' - |X|)/\nu'$, we obtain our result (with $\nu'$ instead of $\nu$).

*Case $\nu > \mu$.* This case is slightly simpler. There exists $x_1$ such that

$$\text{if } |X| \geq x_1, \quad \text{then } \mathbb{E}P(X) \leq \nu.$$

Then, arguing as in the case $\nu < \mu$, we have

$$\mathbb{P}(H(X) \leq n) \leq \mathbb{P}\left(\sum_{k=0}^{n-1} P_k > |X| - x_1\right)$$

$$\leq \mathbb{P}\left(\sum_{k=0}^{n-1} P_k - \mathbb{E}(P_k|\mathcal{F}_k) > |X| - x_1 - \sum_{k=0}^{n-1} \mathbb{E}(P_k|\mathcal{F}_k)\right)$$

$$\leq \mathbb{P}\left(\sum_{k=0}^{n-1} P_k - \mathbb{E}(P_k|X_k) > |X| - n\nu\right)$$

$$\leq C_1 n(|X| - n\nu)^{-r'}. \qquad \square$$

6.3. *Shape of regenerative navigation.* Let $\mathcal{A}$ be a regenerative navigation on $N^O$ and $\theta$ its associated regenerative time. We define $P^\theta(X) = |X| - |X_\theta| = |X| - |\mathcal{A}^\theta(X)|$ and $F_X^\theta$ to be the distribution of $P^\theta(X)$. We assume

(A6.3)
$$\begin{cases} \text{(i)} & \mathcal{A} \text{ is a regenerative navigation with non negative progress;} \\ \text{(ii)} & \sup_{X \in \mathbb{R}^d} \mathbb{E} P^\theta(X)^r < \infty \text{ and } \mathbb{E} \theta^r < \infty \text{ for some } r > d+2; \\ \text{(iii)} & F_X^\theta \text{ converges weakly to a distribution } F^\theta \text{ with } \int r F^\theta(dr) > 0. \end{cases}$$

We let $\overline{\theta} = \lim_{|X| \to \infty} \mathbb{E}\theta(X)$ and $\mu = 1/\overline{\theta} \int r F^\theta(dr) > 0$. From Proposition 3.1 and Lemma 3.7, as $|X|$ tends to infinity, a.s. $H(X)/|X| \to \mu$. Not surprisingly, we obtain the next two results as corollaries of Theorems 6.1 and 6.2.

COROLLARY 6.4. *Under the foregoing assumption* (A6.3), *the conclusions of Theorems 6.2 and 6.1 hold true for* $\mathcal{A}$.

As an example, for the small world navigation, we deduce Theorem 1.8.

PROOF OF COROLLARY 6.4. As in the proof of Theorem 6.1, it is sufficient to show that the conclusions of Theorem 6.2 hold true. Let $(\theta_k)$ denote the regenerative sequence, $\tilde{\mathcal{A}}(X) = \mathcal{A}^\theta(X)$ and $H^\theta(X) = \inf\{k : X_{\theta_k} = 0\} = \sup\{k : \tilde{\mathcal{A}}^k(X) = 0\}$. We first assume that $|X| < n\nu - 1$ and $\nu < \mu$. We may find $0 < \delta < \overline{\theta}$ such that $\nu' = \nu\overline{\theta}/\delta < \mu$ and $|X| < \nu'n - 1$. We obtain

$$\begin{aligned}
(34) \quad \mathbb{P}(H(X) > n) &\leq \mathbb{P}\left(H^\theta(X) > \frac{n}{\delta}\right) + \mathbb{P}(\theta_{\lfloor n/\delta \rfloor} < n) \\
&\leq \mathbb{P}\left(H^\theta(X) > \frac{n}{\delta}\right) + \mathbb{P}\left(\sum_{k=0}^{\lfloor n/\delta \rfloor - 1} \theta_{k+1} - \theta_k < n\right) \\
&\leq \mathbb{P}\left(H^\theta(X) > \frac{n}{\delta}\right) \\
&\quad + \mathbb{P}\left(\left|\sum_{k=0}^{\lfloor n/\delta \rfloor - 1} (\theta_{k+1} - \theta_k - \overline{\theta})\right| > n(1 - \overline{\theta}/\delta) - 1\right).
\end{aligned}$$

Since $|X| < n\nu' - 1 < (\frac{n}{\delta})(\nu'\overline{\theta}) - 1$, Theorem 6.2 applies to $\tilde{\mathcal{A}}$ and $\tilde{\nu} = \nu'\overline{\theta} < \mu\overline{\theta}$. The first term in the latter inequality (34) is thus bounded by $C_1 n \delta^{-1} (n\nu'\overline{\theta}/\delta - |X|)^{-r'} = C_1 n \delta^{-1} (n\nu - |X|)^{-r'}$. We can also apply Lemma 6.3 to the sequence of i.i.d. variables $Y_k = \theta_{k+1} - \theta_k - \overline{\theta}$. Thus, the second term in (34) is upper bounded by $C_1(1 - \overline{\theta}/\delta - 1/n)^{-r'} n^{1-r'}$ for $n$ large enough to guarantee $1 - \overline{\theta}/\delta - 1/n > t_0$ with $0 < t_0 < 1 - \overline{\theta}/\delta$. Finally, we obtain [since $n \geq (n\nu - |X|)/\nu$], for $n$ large enough,

$$\mathbb{P}(H(X) > n) \leq C_1 n(n\nu - |X|)^{-r'}.$$

By suitably increasing $C_1$, we obtain the result for all $n \geq 1$. □

## APPENDIX

**A.1. Proof of Lemma 2.3.** Let $\tilde{V}(X, N^X) = \{Y \in N : \tilde{U}(X,Y) \leq f(|X - Y|)\}$, where the variables $\tilde{U}$ are independent and uniformly distributed on $[0,1]$. We define $\tilde{\rho}(X, N^X) = \inf\{r > 0 : \tilde{V}(X, N^X) \subset B(X, r)\}$. $\tilde{V}(X, N^X) \setminus \{X\}$ is a PPP of intensity $f(|X - x|)\,dx$, hence, for $t \geq 0$,

$$\mathbb{P}(\tilde{\rho}(X, N^X) \geq t) = \mathbb{P}(V(X, N^X) \cap B(X, t)^c = \varnothing) \qquad (35)$$
$$= 1 - e^{-\int_{B(O,t)^c} f(x)\,dx} \leq C_1 t^{d-\beta}.$$

We now define

$$V_O(X, N^{O,X}) = \{Y \in N : U(X, Y) \leq f(|X - Y|)\},$$

where the variables $U$ are independent and uniformly distributed on $[0,1]$, conditioned on the event $\Omega(X) = \{V_O(X, N^{O,X}) \cap B(O, |X|) \neq \varnothing\} = \{P(X) > 0\}$. Similarly, we define

$$V_{-e_1}(X, N^X) = \{X\} \cup \{Y \in N : U_{-e_1}(X, Y) \leq f(|X - Y|)\},$$

where the variables $U_{-e_1}$ are independent and uniformly distributed on $[0,1]$, conditioned on the event $\Omega_{-e_1}(X) = \{V_{e_1}(X, N^X) \cap \mathcal{H}_{-e_1}(X) \neq \varnothing\} = \{P_{-e_1}(X) > 0\}$. We have

$$\mathcal{A}(X) = \arg\inf\{|Y| : Y \in V_O(X, N^{O,X})\}$$

and

$$\mathcal{A}_{-e_1}(X) = \arg\inf\{\langle Y, e_1\rangle : Y \in V_{-e_1}(X, N^X)\}.$$

Let $\rho_{-e_1}(X, N^X) = \inf\{r > 0 : V_{-e_1}(X, N^X) \subset B(X, r)\}$ and $C = \mathbb{P}(\tilde{V}(O, N^O) \cap \mathcal{H}_{-e_1}(O) \neq \varnothing)$. From (35), we deduce

$$\mathbb{P}(\rho_{-e_1}(X, N^X) \geq t) = \frac{\mathbb{P}(\tilde{\rho}(X, N^X) \geq t)}{C} \leq C_1 t^{d-\beta}.$$

Hence, the directed small-world navigation $\mathcal{A}_{-e_1}$ is a.s. defined if and only if $d < \beta$. Similarly, let $\rho_O(X, N^{O,X}) = \inf\{r > 0 : V_O(X, N^{O,X}) \subset B(X, r)\}$ and $C_O(X) = \mathbb{P}(\tilde{V}(X, N^{O,X}) \cap B(O, |X|) \neq \varnothing)$. Since $\inf_{X \in \mathbb{R}^d} C_O(X) > 0$, we obtain, for $\beta > d$,

$$\mathbb{P}(\rho_O(X, N^{O,X}) \geq t) = \frac{\mathbb{P}(\tilde{\rho}(X, N^{O,X}) \geq t)}{C_O(X)} \leq C_1 t^{d-\beta}.$$

We will now couple the variables $U$ and $U_{-e_1}$ via $\tilde{U}$ in order to obtain the statement of the lemma. Without loss of generality, we may assume $\theta > 0$ and $X = xe_x$, with $x > 0$, $e_x \in S^{d-1}$. Let $K(X, e_1) = \{Y \in \mathbb{R}^d : |Y| \geq$ $|X|, \langle Y - X, e_1 \rangle \leq 0\}$, $L(X,e_1) = \{Y \in \mathbb{R}^d : |Y| \leq |X|, \langle Y - X, e_1 \rangle \geq 0\}$ and $M(X,e_1) = B(O,|X|) \setminus L(X,e_1)$. The sets $K(X,e_1)$ and $L(X,e_1)$ are depicted in Figure 2. Set $S_{e_1}(X) = N \cap B(X, \rho_{-e_1}(X, N^X))$ and $S_O(X) = N^O \cap B(X, \rho_O(X, N^{O,X}))$.

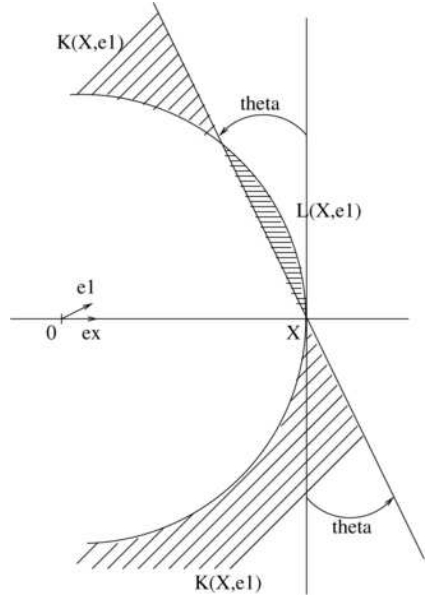

Fig. 2. *The sets $L(X,e_1)$ and $K(X,e_1)$.*

If $S_{e_1}(X) \cap K(X,e_1) = \varnothing$ and $S_O(X) \cap L(X,e_1) = \varnothing$, then both $\mathcal{A}(X)$ and $\mathcal{A}_{-e_1}(X)$ are in $N \cap M(X,e_1)$. Hence, if $S_{e_1}(X) \cap K(X,e_1) = \varnothing$ and $S_O(X) \cap L(X,e_1) = \varnothing$, for all $Y \in N \cap M(X,e_1)$, we may set $U(X,Y) = U_{-e_1}(X,Y) = \tilde{U}(X,Y)$, where $\tilde{U}(X, \cdot)$ is conditioned on the event $\tilde{\Omega}(X,e_1) = \{\tilde{V}(X, N^X) \cap M(X,e_1) \neq \varnothing\}$. It follows that

$$\mathbb{P}(\mathcal{A}_{-e_1}(X) \neq \mathcal{A}(X))$$
(36)
$$\leq \mathbb{P}(S_{e_1}(X) \cap K(X,e_1) \neq \varnothing) + \mathbb{P}(S_O \cap L(X,e_1) \neq \varnothing)$$
$$+ \mathbb{P}(\tilde{\mathcal{A}}_{-e_1}(X) \neq \tilde{\mathcal{A}}(X) | \tilde{\Omega}(X,e_1)),$$

with $\tilde{\mathcal{A}}(X) = \arg\inf\{|Y| : Y \in M(X,e_1) \cap \tilde{V}(X, N^X)\}$ and $\tilde{\mathcal{A}}_{-e_1}(X) = \arg\inf\{\langle Y, e_1 \rangle : Y \in M(X,e_1) \cap \tilde{V}(X, N^X)\}$.

We first upper bound the second term of the right-hand side of (36). Note that $L(X,e_1)$ is contained in a cone of apex $\theta$ (see Figure 2). Letting $\mathcal{C}_\theta$ be a cone issuing from 0 with apex $\theta$, we have

$$\mathbb{P}(S_O(X) \cap L(X,e_1) \neq \varnothing)$$
$$\leq \mathbb{P}(N \cap \mathcal{C}_\theta \cap B(O, \theta^{-1/2}) \neq \varnothing) + \mathbb{P}(\rho_O(X, N^{O,X}) \geq \theta^{-1/2})$$

$$\leq 1 - \exp(-C_0 \theta^{1/2}) + C_1 \theta^{(d-\beta)/2}.$$

The first term of (36) is upper bounded similarly:

$$\mathbb{P}(S_{e_1}(X) \cap K(X, e_1) \neq \varnothing)$$
$$\leq 2\mathbb{P}(N \cap \mathcal{C}_{\arcsin(t/(2|X|))+\theta} \cap B(O, t) \neq \varnothing) + 2\mathbb{P}(\rho_{-e_1}(X, N^X) \geq t)$$
$$\leq 2\left(1 - \exp\left(-C_0 t \left(\arcsin\left(\frac{t}{2|X|}\right) + \theta\right)\right) + C_1 t^{d-\beta}\right).$$

If $1/\sqrt{|X|} \leq \theta$, we choose $t = \sqrt{|X|}$, otherwise we choose $t = 1/\sqrt{\theta}$.

It remains to bound the last term of (36). For $Y \in B(X, \rho(X))$, let $K'(Y, X, e_1) = \{Z \in B(X, \tilde{\rho}(X, N^X)) : |Z| \geq |Y|, \langle Z - Y, e_1 \rangle \leq 0\}$, that is, the set of points with a larger norm but a smaller projection on $e_1$; we have $K'(X, X, e_1) = K(X, e_1)$. We obtain

$$\mathbb{P}(\tilde{\mathcal{A}}_{-e_1}(X) \neq \tilde{\mathcal{A}}(X))$$
$$\leq \mathbb{P}(\exists Y \in N \cap B(X, \tilde{\rho}(X, N^X)) : K'(Y, X, e_1) \cap N \neq \varnothing)$$
$$\leq \mathbb{P}(\tilde{\rho}(X, N^X) \geq t) + \mathbb{P}(N(B(X, t)) \geq n)$$
$$\quad + 2n\mathbb{P}(N \cap \mathcal{C}_{\arcsin(t/(2(|X|-t)))+\theta} \cap B(O, t) \neq \varnothing).$$

We choose $t = \min(|X|, 1/\theta)^{1/3d}$, $n = |X|^{1/2}$. Then, using the inequality $\mathbb{P}(N(B(X,t)) > n) \leq \exp(-n \ln \frac{n}{C_1 t^d})$, we obtain the required bound.

**A.2. Proof of Lemma 2.4.** Since the proof relies on explicit computation and does not involve any subtle arguments, we skip most of the details. Let $\tilde{V}(X, N^X) = \{X\} \cup \{Y \in N : \tilde{U}(X, Y) \leq f(|X - Y|)\}$, where the variables $\tilde{U}$ are independent and uniformly distributed. $\tilde{V}(X, N^X) \setminus \{X\}$ is a nonhomogeneous Poisson point process of intensity $f(|X - x|)\, dx$.

*Statement* 1. Conditioning on the event $\{\tilde{V}(X, N^X) \cap \mathcal{H}_{e_1}(O) \neq \varnothing\}$, we have

$$\mathbb{P}(P_{e_1}(0) > t) = \mathbb{P}(\tilde{V}(O, N^O) \cap \mathcal{H}(te_1) \neq \varnothing \mid \tilde{V}(X, N^X) \cap \mathcal{H}_{e_1}(O) \neq \varnothing)$$
$$= (1 - e^{-\int_{\mathcal{H}_{e_1}(te_1)} f(y)dy})/(1 - e^{-\int_{\mathcal{H}_{e_1}(O)} f(y)dy})$$
$$\sim (1 - e^{-\int_{\mathcal{H}_{e_1}(O)} f(y)dy})^{-1} \int_{\mathcal{H}(te_1)} f(y)\, dy$$

as $t$ tends to infinity. Let $\Lambda_t = \int_{\mathcal{H}(te_1)} f(y)\, dy$. Writing $y = r\cos\theta e_1 + r\sin\theta e_2$, with $\langle e_1, e_2 \rangle = 0$ and $e_2 \in S^{d-1}$, we obtain

$$\Lambda_t = 2\omega_{d-2} \int_0^{\pi/2} \int_{t/\cos\theta}^\infty f(r) r^{d-1}\, dr\, d\theta$$

$$\sim 2\omega_{d-2} \int_0^{\pi/2} \int_{t/\cos\theta}^{\infty} cr^{d-\beta-1} \, dr \, d\theta$$

$$\sim \frac{2c\omega_{d-2}}{\beta - d} t^{d-\beta} \int_0^{\pi/2} \cos^{\beta-d}\theta \, d\theta.$$

*Statement* 2. We can suppose, without loss of generality, that $X = -xe_1$, with $x > 0$. By definition, for $0 \leq t < x$,

$$\mathbb{P}(P(X) > t)$$
$$= \mathbb{P}(\tilde{V}(X, N^{O,X}) \cap B(O, x-t) \neq \varnothing \mid \tilde{V}(X, N^{O,X}) \cap B(O, |X|) \neq \varnothing)$$
$$= \mathbb{P}(\tilde{V}(X, N^{O,X}) \cap B(O, x-t) \neq \varnothing)/\mathbb{P}(\tilde{V}(X, N^{O,X}) \cap B(O, x-0) \neq \varnothing).$$

It thus suffices to compute

$$\mathbb{P}(\tilde{V}(X, N^{O,X}) \cap B(O, x-t) \neq \varnothing)$$
$$= 1 - (1 - f(x)) \exp\left(-\int_{B(O,x-t)} f(|X-y|) \, dy\right).$$

In $\mathbb{R}^2$, for $u \in (0,1)$ and $0 \leq \theta < \arcsin(1-u)$, the straight line with equation $y = \tan\theta$ intersects the sphere of radius $u$ and center $(1,0)$ at two points of respective norms $A(\theta, u)$ and $B(\theta, u)$. A direct computation leads to

$$A(\theta, u) = \cos\theta\left(1 - \sqrt{1 - \frac{u(2-u)}{\cos^2\theta}}\right) = \frac{u}{\cos\theta} + o\left(\frac{u}{\cos\theta}\right)$$

$$B(\theta, u) = \cos\theta\left(1 + \sqrt{1 - \frac{u(2-u)}{\cos^2\theta}}\right) = 2\cos\theta - \frac{u}{\cos\theta} + o\left(\frac{u}{\cos\theta}\right).$$

Let $\Lambda_t(x) = \int_{B(O,x-t)} f(|X-y|) \, dy$. As $t, x$ tend to infinity and $t/x$ tends to 0, we obtain

$$\Lambda_t(x) = 2\omega_{d-2} \int_0^{\arcsin(1-t/x)} \int_{xA(\theta,t/x)}^{xB(\theta,t/x)} f(r) r^{d-1} \, dr \, d\theta$$

$$\sim 2\omega_{d-2} \int_0^{\arcsin(1-t/x)} \int_{xA(\theta,t/x)}^{xB(\theta,t/x)} cr^{d-\beta-1} \, dr \, d\theta$$

$$\sim \frac{2c\omega_{d-2}}{\beta - d} \int_0^{\arcsin(1-t/x)} (xA(\theta, t/x))^{d-\beta} - (xB(\theta, t/x))^{d-\beta} \, d\theta$$

$$\sim \overline{F}(t).$$

It also follows that

$$|\Lambda_t(x) - \Lambda(t)| \leq \int_{\arcsin(1-t/x)}^{\pi/2} f(r) r^{d-1} \, dr \, d\theta$$

$$+ \int_0^{\arcsin(1-t/x)} \int_{t/\cos(\theta)}^{xA(\theta,t/x)} f(r) r^{d-1} \, dr \, d\theta$$

$$+ \int_0^{\arcsin(1-t/x)} \int_{xB(\theta,t/x)}^{\infty} f(r) r^{d-1} \, dr \, d\theta.$$

If $t = x^{d-\beta} \varepsilon(x)$ with $\varepsilon \in \ell^0(\mathbb{R}_+)$, we obtain that $t^{\beta-d} |\Lambda_t(x) - \Lambda(t)|$ tends to 0.

*Statement* 3. Let $Q(X) = |\mathcal{A}(X)|/x^\alpha = (x - P(X))/x^\alpha$, with $|X| = x$ and $\alpha = 1 - (d-\beta)/2 \in (0,1)$. Letting $0 < s < x^{1-\alpha}$, we have

$$\mathbb{P}(Q(X) < s)$$
$$= \mathbb{P}(\tilde{V}(X, N^{O,X}) \cap B(O, sx^\alpha) \neq \varnothing) / \mathbb{P}(\tilde{V}(X, N^{O,X}) \cap B(O, x) \neq \varnothing).$$

Again, it suffices to compute

$$\mathbb{P}(\tilde{V}(X, N^{O,X}) \cap B(O, sx^\alpha) \neq \varnothing)$$
$$= 1 - (1 - f(x)) \exp\left(-\int_{B(O, sx^\alpha)} f(|X-y|) \, dy\right)$$
$$= 1 - (1 - f(x)) \exp(-\Lambda_{x-sx^\alpha}(x)),$$

with, as $x$ tends to $+\infty$, uniformly in $s < x^{1-\alpha'}$, $\alpha' > \alpha$,

$$\Lambda_{x-sx^\alpha}(x) \sim 2\omega_{d-2} \int_0^{\arcsin(sx^{\alpha-1})} \int_{xA(\theta, 1-sx^{\alpha-1})}^{xB(\theta, 1-sx^{\alpha-1})} cr^{d-\beta-1} \, dr \, d\theta$$

$$\sim \frac{2c\omega_{d-2}}{d-\beta} \int_0^{\arcsin(sx^{\alpha-1})} (xB(\theta, 1 - sx^{\alpha-1}))^{d-\beta}$$
$$- (xA(\theta, 1 - sx^{\alpha-1}))^{d-\beta} \, d\theta.$$

We have $B(\theta, 1 - sx^{\alpha-1}) = \cos\theta(1 + \sqrt{s^2 x^{2(\alpha-1)}/\cos^2\theta - \tan^2\theta}) = \cos\theta(1 + \sqrt{s^2 x^{\beta-d}/\cos^2\theta - \tan^2\theta})$ and $A(\theta, 1 - sx^{\alpha-1}) = \cos\theta(1 - \sqrt{s^2 x^{\beta-d}/\cos^2\theta - \tan^2\theta})$. Hence, as $x$ tends to $\infty$, we have

$$(xB(\theta, 1 - sx^{\alpha-1}))^{d-\beta} - (xA(\theta, 1 - sx^{\alpha-1}))^{d-\beta}$$
$$\sim 2(d-\beta) x^{d-\beta} \cos^{d-\beta}\theta \sqrt{s^2 x^{\beta-d}/\cos^2\theta - \tan^2\theta}$$

and we obtain

$$\Lambda_{x-sx^\alpha}(x) \sim 4c\omega_{d-2} \int_0^{\arcsin(sx^{(\beta-d)/2})} x^{d-\beta}$$
$$\times \cos^{d-\beta}\theta \sqrt{s^2 x^{\beta-d}/\cos^2\theta - \tan^2\theta} \, d\theta$$
$$\sim 4c\omega_{d-2} s^2.$$

Finally, we have proven that, uniformly in $s < x^{(d-\beta)/2-\eta}$ (for some $\eta > 0$),

$$\lim_{|X|\to\infty} \mathbb{P}(Q(X) > s) = \exp(-4c\omega_{d-2}s^2). \tag{37}$$

*Statement* 4. Again, we suppose that $X = -xe_1$, with $x > 0$. Letting $s > 0$ and $u = 1 - \exp(-s) \in (0,1)$, we have

$$\mathbb{P}(\tilde{P}(X) > s)$$
$$= \mathbb{P}(P(X) > xu)$$
$$= (1 - (1-f(x))e^{-\int_{B(O,(1-u)x)} f(|X-y|)dy})/\mathbb{P}(\tilde{V}(X, N^{O,X}) \cap B(O,x) \neq \emptyset).$$

As above, with $\Lambda_t(x) = \int_{B(O,x-t)} f(|X-y|)\,dy$, we have

$$\Lambda_{ux}(x) = 2\omega_{d-2} \int_0^{\arcsin(1-u)} \int_{xA(\theta,u)}^{xB(\theta,u)} f(r) r^{d-1}\,dr\,d\theta$$
$$\sim 2\omega_{d-2} \int_0^{\arcsin(1-u)} \int_{xA(\theta,u)}^{xB(\theta,u)} c/r\,dr\,d\theta$$
$$\sim 2c\omega_{d-2} \int_0^{\arcsin(1-u)} \ln\frac{B(\theta,u)}{A(\theta,u)}\,d\theta.$$

We define

$$\tilde{F}(s) = 1 - \exp\left(-2c\omega_{d-2} \int_0^{\arcsin(\exp(-s))} \ln\frac{B(\theta, 1-\exp(-s))}{A(\theta, 1-\exp(-s))}\,d\theta\right) \tag{38}$$
$$= 1 - \exp\left(-c\int_{B(O,\exp(-s))} |e_1 - y|^{-d}\,dy\right)$$
$$\sim 4c\omega_{d-2}e^{-2s} \tag{39}$$

as $s$ tends to $+\infty$. Statement 4 follows.

**A.3. Tail inequality in the GI/GI/$\infty$.** Let $\{\sigma_n, \tau_n\}, n \in \mathbb{Z}$ be an i.i.d. sequence of $\mathbb{R}_+ \times \mathbb{R}_+$-valued random variables representing the service times and interarrival times in a GI/GI/$\infty$ queue. The random variables $(\sigma_n)$ and $(\tau_n)$ are independent. We set $T_0 = 0$ as the arrival time of customer 0; for $n \geq 1$, $T_n = \sum_{k=0}^{n-1} \tau_k$ is the arrival time of the $n$th customer. Let $Y \in \mathbb{R}_+$ be a nonnegative initial condition, independent of the $\{\sigma_n, \tau_n\}$ sequence. We set $W_0^{[Y]} = Y$ and, for $n \geq 1$, we define

$$W_n^{[Y]} = \max\left(Y - (T_n - T_0), \max_{1 \leq i \leq n} \sigma_{i-1} - (T_n - T_i)\right)^+$$
$$= \max(W_{n-1}^{[Y]} - \tau_{n-1}, \sigma_{n-1})$$

NAVIGATION ON A POISSON POINT PROCESS 39

(where, by convention, $\sum_{k=n}^{n-1} \cdot = 0$). The random variable $W_n^{[Y]}$ is the largest residual service time just after the arrival of the $n$th customer in the GI/GI/$\infty$ queue with initial condition $Y$. Let $\mathcal{F}_n$ be the $\sigma$-field generated by the random variables $Y$ and $\{(\sigma_k, \tau_k), k = 0, \ldots, n-1\}$. Consider the $\{\mathcal{F}_n\}$-stopping time

$$(40) \qquad \theta(Y) = \theta = \inf\{n \geq 1 : W_n^{[Y]} = 0\}.$$

$\theta$ is the time needed to empty all queues. The proof of the next lemma follows from a classical computation in queueing theory; see [7].

LEMMA A.1. *Let $\theta$ be the stopping time defined in (40). Assume:*

(i) *there exists an $\alpha > 1$ such that $\mathbb{P}(\sigma_1 > t) \leq C_1 t^{-\alpha}$, $\mathbb{P}(Y > t) \leq C_1 t^{1-\alpha}$;*
(ii) $\mathbb{P}(\tau_1 > 0) > 0$;
(iii) $\mathbb{P}(\sigma_1 = 0) > 0$.

*$\theta$ is then a.s. finite and if $\alpha > 2$, there exists $C_1 > 0$ such that*

$$\mathbb{E}\theta < \infty \quad \text{and} \quad \mathbb{P}(\theta > t) \leq C_1 t^{2-\alpha}.$$

*Moreover, for all $n$, $W_n^{[0]} \leq_{st} M$, where $M$ is the stationary workload and $\mathbb{P}(M > t) \leq C_1 t^{1-\alpha}$.*

*If assumption (i) is replaced by*

(i′) *there exists $s > 0$ such that $\mathbb{E}\exp(s\sigma_1) < \infty$ and $\mathbb{E}\exp(sY) < \infty$,*

*then there exists $C_0 > 0$ such that $\mathbb{E}\exp(C_0\theta) < \infty$ and $\mathbb{E}\exp(C_0 M) < \infty$.*

**Acknowledgments.** The author thanks Neil O'Connell and the Boole Center. This work has benefited from enlightening discussions with François Baccelli.

DEPARTMENTS OF EECS AND STATISTICS
UNIVERSITY OF CALIFORNIA
BERKELEY, CALIFORNIA 94720-1770
USA
E-MAIL: [charles.bordenave@ens.fr](charles.bordenave@ens.fr)